\documentclass[11pt]{amsart}

\usepackage{epsfig}
\usepackage[]{rotating}

\newtheorem{thm}{Theorem}[section]
\newtheorem{cor}[thm]{Corollary}
\newtheorem{lem}[thm]{Lemma}

\theoremstyle{remark}
\newtheorem{rem}[thm]{Remark}
\newtheorem{defin}[thm]{Definition}
\newtheorem{exam}[thm]{Example}

\newtheorem{nota}[thm]{Notation}

\newtheorem{prp}[thm]{Proposition}
\newtheorem{dfn}[thm]{Definition}

\newcommand{\N}{{\mathbb{N}}}

\newcommand{\grad}{{\textrm{grad }}}
\newcommand{\ric}{{\textrm{Ric}}}
\newcommand{\bcwp}{$(\psi,\mu)$-\emph{bcwp}}
\newcommand{\bcwpp}{\emph{bcwp }}
\newcommand{\bcwpar}[1]{$(#1)$-\emph{bcwp}}

\newcommand{\discr}[1]{\emph{discr }(#1)}
\newcommand{\twistpar}[4]{$#1 \times_{(#3;#4)} #2$}

\numberwithin{equation}{section}
\renewcommand{\theequation}{\arabic{section}.\arabic{equation}}

\begin{document}

\title{CURVATURE IN SPECIAL BASE CONFORMAL WARPED PRODUCTS}

\author{Fernando Dobarro \\
\& \\
B\"{u}lent \"{U}nal}

\address[F. Dobarro]{Dipartimento di Matematica e Informatica,
Universit\`{a} degli Studi di Trieste, Via Valerio 12/b, I-34127
Trieste, Italy} \email {dobarro@dmi.units.it}

\address[B. \"{U}nal]{Department of Mathematics, Bilkent University,
         Bilkent, 06800 Ankara, Turkey}
\email{bulentunal@mail.com}



\date{April 25, 2008}



%
%


\begin{abstract}
We introduce the concept of a base conformal warped product of two
pseudo-Riemannian manifolds. We also define a subclass of this
structure called as a special base conformal warped product.
%
%
After, we explicitly mention many of the relevant fields where
metrics of these forms and also considerations about their
curvature related properties play important rolls. Among others,
we cite general relativity, extra-dimension, string and
super-gravity theories as physical subjects and also the study of
the spectrum of Laplace-Beltrami operators on p-forms in global
analysis.
%
%
Then, we give expressions for the Ricci tensor and scalar
curvature of a base conformal warped product in terms of Ricci
tensors and scalar curvatures of its base and fiber, respectively.
Furthermore, we introduce specific identities verified by
particular families of, either scalar or tensorial, nonlinear
differential operators on pseudo-Riemannian manifolds. The latter
allow us to obtain new interesting expressions for the Ricci
tensor and scalar curvature of a special base conformal warped
product and it turns out that not only the expressions but also
the analytical approach used are interesting from the physical,
geometrical and analytical point of view.
%
%
Finally, we analyze, investigate and characterize possible
solutions for the conformal and warping factors of a special base
conformal warped product, which guarantee that the corresponding
product
is Einstein. Besides all, we apply these results to a generalization
of the Schwarzschild metric.
\end{abstract}

\maketitle

 \pagebreak


\centerline{\textsc{Contents}}

\begin{list}{\item contents}{}
  \item[\ref{sec:Introduction}.]Introduction
  \item[\ref{sec:Some differential operators}.]Some families of
differential operators
  \item[\ref{sec:base conformal warped products}.]About base
conformal warped products
   \begin{itemize}
  \item[\ref{subsec:Covariant Derivatives}.]Covariant Derivatives
  \item[\ref{subsec:Riemannian Curvatures}.]Riemannian Curvatures
  \item[\ref{subsec:Ricci Curvatures}.]Ricci Curvatures
  \item[\ref{subsec:Scalar Curvature}.]Scalar Curvature
   \end{itemize}
  \item[\ref{sec:Curvature of sbcwp}.]Curvature of $({(B \times
F)}_{m+k},\psi^{2\mu} g_{B} + \psi^{2}g_{F})$
\begin{itemize}
  \item[\ref{subsection sbcwp Ricci}.] Ricci Tensor
  \item[\ref{subsec: sbcwp scalar curvature}.] Scalar curvature
\end{itemize}
  \item[\ref{sec:constant scalar curvature bcwp with m ge 2}.]The
nonlinearities in the \bcwp  $\,$ scalar curvature relations
   \begin{itemize}
     \item[\ref{subsec: General Properties}.]Base $B_m$ with
dimension $m \ge 2$
     \item[\ref{subsec:constant scalar curvature bcwp with m =
1}.]Base $B_m$ with dimension $m = 1$
\end{itemize}
  \item[\ref{sec:Some Examples}.]Some Examples and Remarks
  \item[\ref{sec:Conclusions}.]Conclusions and Future Directions
  \item[\ref{sec:appendix}.]Appendix
  \item[$\qquad\:$References]
\end{list}

\bigskip




\section{Introduction}
\label{sec:Introduction}

The main concern of the present paper is so called base conformal
warped products (for brevity, we call a product of this class as a
{\it bcwp}) and their interesting curvature related geometric
properties. One can consider {\it bcwp}'s as a generalization of
the classical singly warped products. Before we mention physical
motivations and applications of
{\it bcwp}'s, we will explicitly define warped products and
briefly mention their different types of extensions.
This is the first of a series of articles where we deal with the
study of curvature questions in {\it bcwp}'s, the latter also give
rise to interesting problems in nonlinear analysis.
%

\bigskip

Let $B=(B_m,g_{B})$ and $F=(F_k,g_{F})$ be two pseudo-Riemannian
manifolds of dimensions $m \geq 1$ and $k \geq 0,$ respectively and
also let $B \times F$ be the usual product manifold of $B$ and $F$.
Given a smooth function $w \in C^\infty_{>0}(B)=\{v \in C^\infty(B):
v>0\}$, the \emph{warped product} $B \times _w F = ({(B \times_w
F)}_{m+k},g=g_{B} + w^2g_{F})$ was first defined by Bishop and
O'Neill in \cite{BishopONeil69} in order to study manifolds of
negative curvature. Moreover, they obtained expressions for the
sectional, Ricci and scalar curvatures of a warped product in terms
of sectional, Ricci and scalar curvatures of its base and fiber,
respectively (see also \cite{Beem-Ehrlich81,Beem-Ehrlich-Easley96,
Beem-Ehrlich-Powell82, Besse87,ONeil83} and for other developments
about warped products see for instance \cite{Chen02,
DajczerTojeiro04,DostoglouEhrlich04,Mirzoyan03,Zeghib99a, Zeghib99b,
Zeghib01}).

\bigskip

From now on, we will use the Einstein summation convention over
repeated indices and consider only connected manifolds.
Furthermore, we will denote the Laplace-Beltrami operator on
$(B,g_{B})$ by $\Delta_{B}(\cdot),$ i.e.,

$$\Delta_{B}(\cdot)={\nabla^{B}}^i{\nabla^{B}}_i (\cdot)=
\displaystyle{\frac{1}{\sqrt{\mid g_{B} \mid}}} {\partial_i}
\big({\sqrt{\mid g_{B} \mid}}g_{B}^{ij}{\partial_j}(\cdot)\big).$$

\bigskip

Note that $\Delta_{B}$ is elliptic if $(B,g_{B})$ is Riemannian
and it is hyperbolic when $(B,g_{B})$ is Lorentzian. If
$(B,g_{B})$ is neither Riemannian nor Lorentzian, then the
operator is called as ultra-hyperbolic (see \cite{Branson97}).

\bigskip

In \cite{PongeReckziegel93}, Ponge and Reckziegel generalized the
notion of warped product to twisted and doubly-twisted products,
i.e., a \emph{doubly-twisted product} $B \times_{(\psi_{0};
\psi_{1})} F$ can be defined as the usual product $B \times F$
equipped with the pseudo-Riemannian metric $\psi_{0}^{2} g_{B} +
\psi_{1}^2 g_{F}$ where $\psi_{0}, \psi_{1} \in C^\infty_{>0}(B
\times F).$ In the case of $\psi_{0} \equiv 1,$ the corresponding
\emph{doubly-twisted product} is called as a \emph{twisted
product} by  B.-Y. Chen (see \cite{Bishop72,Chen81}). Clearly, if
$\psi_{1}$ only depends on the points of $B$, then $B \times_{(1;
\psi_{1})} F$ becomes a warped product. One can also find other
interesting generalizations in \cite{DobarroUnal04,
KatanaevKloschKummer99, Unal93,Unal00,UnalPHD}.

\bigskip

We recall that a pseudo-Riemannian manifold $(B_m,g_{B})$ is
conformal to the pseudo-Riemannian manifold $(B_m,\tilde{g}_{B})$,
if and only if there exists $\eta \in C^\infty (B)$ such that
$\tilde{g}_{B} = e^{\eta} g_{B}$.

\bigskip

From now on, we will call a doubly twisted product as a \emph{base
conformal warped product}  when the functions $\psi_{0}$ and
$\psi_{1}$ only depend on the points of $B.$ For a precise
definition, see \S \ref{sec:base conformal warped products}. In
this article, we deal with \emph{bcwp}'s, and especially with a
subclass called as {\it special base conformal warped products},
briefly \emph{sbcwp}, which can be thought as a mixed structure of
a conformal change in the metric of the base and a warped product,
where there is a specific type of relation between the conformal
factor and the warping function. Precisely, a {\it special base
conformal warped product} is the usual product manifold $B_{m}
\times F_{k}$ equipped with pseudo-Riemannian metric of the form
$\psi^{2\mu}g_{B}+ \psi^{2}g_{F}$ where $\psi \in
C^\infty_{>0}(B)$ and a parameter $\mu \in \mathbb{R}.$ In this
case, the corresponding {\it special base conformal warped
product} is denoted by \bcwp. Note that when $\mu = 0,$ we have a
usual warped product and when $k=0$ we have a usual conformal
change in the base (the fiber is reduced to a point) and if $\mu =
1$ we are in the presence of a conformal change in the metric of a
usual product pseudo-Riemannian manifold.

\bigskip

We remark here that a \emph{sbcwp} can be expressed as a special
conformal metric in a particular warped product, i.e.
$$\displaystyle \psi_{0}^{2} g_{B} + \psi_{1}^2 g_{F} =
\psi_{0}^{2}\left( g_{B} + \frac{\psi_{1}^2}{\psi_{0}^2}
g_{F}\right),$$ where $\psi_{0},\psi_{1}\in C^\infty_{>0}(B)$.

\bigskip

Metrics of this type have many applications in several topics from
the areas of differential geometry, cosmology, relativity, string
theory, quantum-gravity, etc. Now, we want to mention some of the
major ones.

\begin{description}
\item[i] In the construction of a large class of non trivial
static anti de Sitter vacuum space-times
\begin{itemize}
  \item
In the Schwarzschild solutions of the Einstein equations
\begin{equation}\label{eq:Schwarzschild-3}
  ds^{2}=-\left(1-\frac{2M}{r}\right)dt^{2} +
  \frac{1}{\displaystyle
1-\frac{2M}{r}}dr^{2}+r^{2}(d\theta^{2}+\sin ^{2} \theta
d\phi^{2})
\end{equation}
(see \cite{Anderson-Chrusciel-Delay02,Besse87,
Hawking-Ellis73,ONeil83,Schwarzschild1916,Thorne03}).
  \item
In the Riemannian Schwarzschild metric, namely
\begin{equation}\label{eq:Schwarzschild-1a}
  (\mathbb{R}^{2} \times \mathbf{S}^{2}, g_{Schw}),
\end{equation}
where
\begin{equation}\label{eq:Schwarzschild-1b}
 g_{Schw}= u^{2}d\phi^{2}+u^{-2}dr^{2}+r^{2}g_{S^{2}(1)}
\end{equation}
and $\displaystyle u^{2}= 1 + r^{2} -\frac{2m}{r}$, $m>0$ (see
\cite{Anderson-Chrusciel-Delay02}).
  \item
In the \emph{``generalized Riemannian anti de Sitter
$\mathbf{T}^{2}$ black hole metrics"} (see \S3.2 of
\cite{Anderson-Chrusciel-Delay02} for details).
\end{itemize}

\bigskip

Indeed, let $(F_2,g_F)$ be a pseudo-Riemannian manifold and $g$ be
a pseudo-metric on $\mathbb{R}_+ \times \mathbb{R} \times F_2$
defined by
\begin{equation}\label{eq:Schwarzschild to bcwp-1}
g= \frac{1}{u^{2}(r)}dr^{2} \pm u^{2}(r) dt^2 + r^{2}g_F.
\end{equation}
After the change of variables $s=r^2$, $\displaystyle
y=\frac{1}{2} t$ and hence $ds^2=4 r^2 dr^2$ and $ \displaystyle
dy^2 = \frac{1}{4} dt^2$. Then \eqref{eq:Schwarzschild to bcwp-1}
is equivalent to
\begin{equation}
\label{eq:Schwarzschild to bcwp-2}
\begin{split}
g&= \frac{1}{\sqrt{s}}\Big[\frac{1}{4 \sqrt{s}
u^{2}(\sqrt{s})}ds^{2} \pm 4 \sqrt{s} u^{2}(\sqrt{s}) dy^2\Big] +
s g_F \\
&= (s^\frac{1}{2})^{2(-\frac{1}{2})}
\Big[(2 s^\frac{1}{4}  u(s^\frac{1}{2}))^{2(-1)}ds^{2} \pm
(2 s^\frac{1}{4}  u(s^\frac{1}{2}))^2 dy^2\Big]
 + (s^\frac{1}{2})^{2} g_F.
\end{split}
\end{equation}
Note that roughly speaking, $g$ is a nested application of two
\bcwp's. That is, on $\mathbb{R}_+ \times \mathbb{R}$ and taking
\begin{equation}\label{eq:Schwarzschild to bcwp-3}
  \psi_1(s) = 2 s^\frac{1}{4}  u(s^\frac{1}{2}) \textrm{ and }
\mu_1=-1,
\end{equation}
the metric inside the brackets in the last member of
\eqref{eq:Schwarzschild to bcwp-2} is a \bcwpar{\psi_1,\mu_1},
while the metric $g$ on $(\mathbb{R}_+ \times \mathbb{R}) \times
F_2$ is a \bcwpar{\psi_2,\mu_2} with
\begin{equation}\label{eq:Schwarzschild to bcwp-4}
  \psi_2(s,y)=s^\frac{1}{2} \textrm{ and } \mu_2 = -\frac{1}{2}.
\end{equation}

\item[ii] In the Ba\~{n}ados-Teitelboim-Zanelli (BTZ) and de
Sitter (dS) black holes (see
\cite{Aharony-Gubser-Maldacena-Ooguri-Oz,
Banados-Teitelboim-Zanelli92,
Banados-Henneaux-Teitelboim-Zanelli93,
DobarroUnal04,Hong-Choi-Park03,Peet2000} for details).

\item[iii] In the study of the spectrum of Laplace-Beltrami
operator for $p-$forms. For instance in Equation $(1.1)$ of
\cite{Antoci03}, the author considers the structure that follows:
let $\overline{M}$ be an $n$-dimensional compact, Riemannian
manifold with boundary, and let $y$ be a boundary-defining
function; she endows the interior $M$ of $\overline{M}$ with a
Riemannian metric $ds^2$ such that in a small tubular neighborhood
of $\partial M$ in $M$, $ds^2$ takes the form
\begin{equation}\label{eq:antoci}
ds^2 = e^{-2(a+1)t} dt^2 + e^{-2bt}d\theta^2_{\partial M},
\end{equation}
where $t := -\log y \in (c,+\infty)$ and $d\theta^2_{\partial M}$
is the Riemannian metric on $\partial M$ (see \cite{Antoci03,
Melrose95} and references therein for details).

\item[iv] In the Kaluza-Klein theory (see \cite[\S 7.6,
Particle Physics and Geometry]{Wesson99} and
\cite{OverduinWesson98}) and in the Randall-Sundrum theory
\cite{Frolov01,Greene-Schalm-Shiu00,Randall-Sundrum99a,Randall-Sundrum99b,Randjbar
Daemi-Rubakov04,Soda02} with $\mu$ as a free parameter. For example
in \cite{Ito01} the following metric is considered
\begin{equation}\label{eq:Ito1}
  e^{2\mathcal{A}(y)} g_{ij}dx^{i}dx^{j} + e^{2\mathcal{B}(y)}dy^{2},
\end{equation}
with the notation $\{x^{i}\}$, $i= 0,1,2,3$ for the coordinates in
the 4-dimensional space-time and $x^{5}=y$ for the fifth
coordinate on an extra dimension. In particular, Ito takes the
ansatz
\begin{equation}\label{eq:Ito2}
  \mathcal{B}=\alpha \mathcal{A},
\end{equation}
which corresponds exactly to our \emph{sbcwp} metrics,
considering $g_{B}=dy^{2}$, $g_{F}=g_{ij}dx^{i}dx^{j}$,
$\displaystyle \psi(y) = e^{\frac{\mathcal{B}(y)}{\alpha}} =
e^{\mathcal{A}(y)}$ and $\displaystyle \mu = \alpha$.

\item [v] In String and Supergravity theories, for instance, in
the Maldacena conjecture about the duality between
compactifications of M/string theory on various Anti-de Sitter
space-times and various conformal field theories (see
\cite{Maldacena98,Maldacena99,Petersen99}) and in warped
compactifications (see \cite{Greene-Schalm-Shiu00, Strominger86}
and references therein). Besides these, there are also frequent
occurrences of this type of metrics in string topics (see \cite{
Gauntlett-Kim-Waldram01, Gauntlett-Kim-Waldram01-b,
Gauntlett-Kim-Pakis-Waldram02,
Gauntlett-Martelli-Sparks-Waldram04,GhezelbashMann04, Lidsey00,
PapadopoulosTownsend99,Soda02} and also
\cite{Aharony-Gubser-Maldacena-Ooguri-Oz, Argurio98, Peet2000,
Schmidt04} for some reviews about these topics).

\item [vi] In the discussion of Birkhoff-type theorems (generally
speaking these are the theorems in which the gravitational vacuum
solutions admit more symmetry than the inserted metric ansatz,
(see \cite[page 372]{Hawking-Ellis73} and  \cite[Chapter
3]{Beem-Ehrlich-Easley96}) for rigorous statements), especially in
Equation 6.1 of \cite{Schmidt97} where, H-J. Schmidt considers a
special form of a \emph{bcwp} and basically shows that if a
\emph{bcwp} of this form is Einstein, then it admits one Killing
vector more than the fiber has. In order to achieve that, the
author considers for a specific value of $\mu$, namely $\mu =
(1-k)/2$, the following problem:

\emph{Does there exist a smooth function $\psi \in
C^\infty_{>0}(B)$ such that the corresponding \bcwpar{\psi,\mu}
$(B_2 \times F_k,\psi^{2 \mu}g_{B} + \psi^{2} g_{F})$ is an
Einstein manifold?} (see also \textbf{(Pb-Eins.)} below.)

\item [vii] In questions of equivariant isometric embeddings
(see \cite{Giblin-Hwang04}).

\item [viii] In the study of bi-conformal transformations,
bi-conformal vector fields and their applications (see
\cite[Remark in Section 7]{GarciaParrado-Senovilla04} and
\cite[Sections 7 and 8]{GarciaParrado04}).

\end{description}

\bigskip

In order to study the curvature of \bcwp's we organized the paper
as follows:

\bigskip

In \textbf{\S  \ref{sec:Some differential operators}}, we study a
specific type of homogeneous non-linear second order partial
differential operator closely related to those  with terms
including $\|\nabla^{B} (\cdot)\|_{B}^{2}=g_{B}(\nabla^{B}(\cdot),
\nabla^{B}(\cdot))$ and a generalization where the Hessian tensor
is involved. Operators with this structure are frequent in
physics, differential geometry and analysis (see
\cite{Aubin82,Aubin98,DobarroLamiDozo87,
DobarroUnal04,Hameiri83,Lelong-Ferrand76,
Yoshida-Mahajan-Ohsaki04}).

\bigskip

In \textbf{\S \ref{sec:base conformal warped products}}, we define
precisely the base conformal warped products, compute their
covariant derivatives and Riemann curvature tensor, Ricci tensor
and scalar curvature.

\bigskip

In \textbf{\S \ref{sec:Curvature of sbcwp}}, applying the results
of \textbf{\S \ref{sec:Some differential operators}} we find a
useful formula for the relation among the Ricci tensors
(respectively the scalar curvatures) in a \bcwp. The principal
results of this section are \textbf{Theorem \ref{thm:ricci sbcwp m
ge 3}}, about the Ricci tensor, and the theorem that follows about
the scalar curvature.

\bigskip

\begin{thm} \label{thm:scurv conf warped m ge 2}
Let $B = (B_m,g_{B})$ and $F = (F_k,g_{F})$ be two
pseudo-Rieman\-nian manifolds with dimensions $m \geq 2$ and $k
\geq 0,$ respectively. Suppose that $S_B$ and $S_F$ denote the
scalar curvatures of $B = (B_m,g_{B})$ and $F = (F_k,g_{F}),$
respectively. If $\displaystyle \mu \in \mathbb{R}$ is a parameter
and $\psi \in C^\infty_{>0}(B)$ is a smooth function then, the
scalar curvature $S$ of the base conformal warped product $(B
\times F,g = \psi^{2\mu} g_{B} + \psi^{2}g_{F}) $ verifies,
\begin{description}
  \item[i] If $\displaystyle \mu \neq -\frac{k}{m-1}$, then
\begin{equation}
-\beta \Delta_{B}u + S_{B} u = S u^{2 \mu \alpha + 1}   -
S_{F}u^{2(\mu - 1)\alpha + 1} \label{eq:sc principal 1}
\end{equation}
where
\begin{equation}
\label{eq:alpha principal}
 \alpha =
\frac{2[k+(m-1)\mu]}{\{[k+(m-1)\mu]+(1-\mu)\}k + (m-2)\mu
[k+(m-1)\mu]},
\end{equation}
\begin{equation}
\label{eq:beta principal}
  \beta = \alpha 2[k+(m-1)\mu]>0
\end{equation}
and $\psi = u^{\alpha}>0$.
   \item[ii] If $\displaystyle \mu = -\frac{k}{m-1}$, then
\begin{equation}
\label{eq:sc special}
  -\left[-k^{2} \frac{m-2}{m-1} + k(k+1)\right] \frac{|\nabla^{B}
  \psi|_{B}^{2}}{\psi^{2}}= \psi^{-2\frac{k}{m-1}} S - S_{B} -
    S_{F} \psi^{-2\left(\frac{k}{m-1} + 1\right)}.
\end{equation}
\end{description}

For the case of $m=1$ see Remark \ref{rem:mu special}.

\end{thm}

\bigskip

The relation among the scalar curvatures in a warped product $B
\times _w F$ is given by
\begin{equation}
S=-2k \frac{\Delta_{B}w}{w} -
k(k-1)\frac{g_{B}(\nabla^{B}w,\nabla^{B}w) }{w^{2}} + S_{B} +
\frac{S_{F}}{w^{2}},
    \label{eq:scalar curv wp-1}
\end{equation}
where $\Delta_{B} $ is the Laplace-Beltrami operator on
$(B,g_{B})$ and $S_{B}$, $S_{F}$ and $S$ are the scalar curvatures
of $B$, $F$ and $B \times_w F$, respectively.

In the articles \cite{Dobarro87,DobarroLamiDozo87} the authors
transformed equation \eqref{eq:scalar curv wp-1} into
\begin{equation}
    - \displaystyle\frac{4k}{k+1} \Delta_{B} u + S_{B} u +
S_{F}u^{1-\frac{4}{k+1}}
    = S u,
    \label{eq:scalar curv wp}
\end{equation}
where  $w = u^{{2}\over{k+1}}$ and $u \in C^\infty_{>0}(B)$. Note
that this result corresponds to the case of $\mu=0$ in
\textbf{Theorem \ref{thm:scurv conf warped m ge 2}}.

\bigskip

On the other hand, under a conformal change on the metric of a
pseudo-\\Riemannian manifold $B = (B_m,g_{B})$, i.e.,
$\tilde{g}_{B} = e^{\eta} g_{B}$ with $\eta \in C^\infty(B)$, the
scalar curvature $\tilde{S}_{B}$ associated to the metric
$\tilde{g}_{B}$ is related with the scalar curvature $S_{B}$ by
the equation
\begin{equation}
e^{\eta} \tilde{S}_{B} = S_{B} -(m-1) \Delta_{B}\eta -
(m-1)\frac{m-2}{4} g_{B}(\nabla^{B}\eta,\nabla^{B}\eta).
    \label{eq:yamabe-1}
\end{equation}
When $m \geq 3$, the previous equation becomes
\begin{equation}
    - \displaystyle 4 \frac{m-1}{m-2} \Delta_{B} \varphi + S_{B}
\varphi
    = \tilde{S}_{B} \varphi^{1+\frac{4}{m-2}},
    \label{eq:yamabe}
\end{equation}
where $\tilde{g}_{B}= \varphi ^{\frac{4}{m-2}} g_{B}$ and $\varphi
\in C^\infty_{>0}(B)$.

\bigskip

There is an extensive number of publications about equation
\eqref{eq:yamabe} (see \cite{Aubin82, Aubin98, Branson97,
Caffarelli-Gidas-Spruck89, Chen-Lin97, Chen-Lin01, Escobar92,
Escobar92-b, Guilfoyle-Nolan98, Hebey99, Hebey99-b, Hebey03,
Kazdan83, Kazdan85, Lee-Parker87, Schoen84, Schoen-Yau88}),
especially due to its close relation with the so called Yamabe
problem (see the original Yamabe's article \cite{Yamabe63} and the
related questions posed by Tr\"{u}dinger \cite{Trudinger68}), namely

\medskip
\begin{center}
\begin{minipage}{10 cm}
\label{(Ya):problem}
    \item[\textbf{(Ya)}]\cite{Yamabe63}
    Does there exist a smooth function
    $\varphi \in C^\infty_{>0}(B)$ such that
    $(B,\varphi^{\frac{4}{m-2}}g_{B})$ has
    constant scalar curvature?
\end{minipage}
\end{center}
\medskip

Analogously, in several articles the following problem has been
studied (see \cite{Badiale-Dobarro96,
Chabrowski-doO02,Dobarro87,DobarroLamiDozo87,
DobarroUnal03,Duggal02, Ehrlich-Jung-Kim96,
Ehrlich-Jung-Kim-Shin01,Leung98,Yang98} among others).

\medskip
\begin{center}
\begin{minipage}{10 cm}
\label{(cscwp):problem}
    \item[\textbf{(cscwp)}]
    Is there a smooth function $w \in C^\infty_{>0}(B)$
    such that the warped product
    $B \times_{w} F$ (or equivalently $B \times_{(1,
    w)} F$) has constant scalar curvature?
\end{minipage}
\end{center}
\medskip

The Yamabe problem needs the study of the existence of positive
solutions of equation \eqref{eq:yamabe} with a constant $\lambda
\in \mathbb{R}$ instead of $\tilde{S}_{B}$. On the other hand, the
constant scalar curvature problem in warped products brings to the
study of the existence of positive solutions of the equation
\eqref{eq:scalar curv wp} with a parameter $\lambda \in
\mathbb{R}$ instead of $S$.

\noindent Inspired by these, we propose a mixed problem between
\textbf{(Ya)} and \textbf{(cscwp)}, namely:

\medskip
\begin{center}
\begin{minipage}{10 cm}
  \item[\textbf{(Pb-sc)}] Given $\mu \in \mathbb{R}$, does there exist
 $\psi \in C^\infty_{>0}(B)$ such that the \bcwpar{\psi,\mu}
 $({(B \times F)}_{m+k},\psi^{2 \mu}g_{B} +
 \psi^{2} g_{F})$ has constant scalar curvature?
\end{minipage}
\end{center}
\medskip


\bigskip

\textsl{Note that when $\mu = 0,$ \textbf{(Pb-sc)} corresponds to
the problem \textbf{(cscwp)}, whereas when the dimension of the
fiber $k=0$ and $\mu = 1$, then \textbf{(Pb-sc)} corresponds to
\textbf{(Ya)} for the base manifold. Finally \textbf{(Pb-sc)}
corresponds to \textbf{(Ya)} for the usual product metric with a
conformal factor in $C^\infty_{>0}(B)$ when $\mu=1$.}

\bigskip

Under the hypothesis of \textbf{Theorem \ref{thm:scurv conf warped
m ge 2} i}, the analysis of the problem \textbf{(Pb-sc)} brings to
the study of the existence and multiplicity of positive solutions
$u \in C^\infty_{>0}(B)$ of
\begin{equation}
   -\beta \Delta_{B}u + S_{B} u =
  \lambda u^{2 \mu \alpha + 1} - S_{F}u^{2(\mu - 1)\alpha + 1},
\label{eq:scalar curv bcwp1}
\end{equation}
where all the components of the equation are like in \textbf{Theorem
\ref{thm:scurv conf warped m ge 2} i} and $\lambda$ (the conjectured
constant scalar curvature of the corresponding \emph{sbcwp}) is a
real parameter. We observe that an easy argument of separation of
variables, like in \cite[Section 2]{CotiZelatiDobarroMusina97} and
\cite{DobarroLamiDozo87}, shows that there exists a positive
solution of \eqref{eq:scalar curv bcwp1} only if the scalar
curvature of the fiber $S_H$ is constant. Thus this will be a
natural assumption in the study of \textbf{(Pb-sc)}.

\noindent Furthermore note that the involved nonlinearities in the
right hand side of \eqref{eq:scalar curv bcwp1} dramatically
change with the choice of the parameters, the analysis of these
changes is the subject of \S\ref{sec:constant scalar curvature
bcwp with m ge 2}.

\noindent By taking into account the above considerations and the
scalar curvature results obtained in this article, we will
consider the study of \textbf{(Pb-sc)} and in particular, the
questions mentioned above which are related to the existence and
multiplicity of solution of \eqref{eq:scalar curv bcwp1} in our
forthcoming articles (see \cite{DobarroUnal07}). Let us mention
here that there are several partial results about semi-linear
elliptic equations like \eqref{eq:scalar curv bcwp1} with
different boundary conditions, see for instance
\cite{Alama99,Ambrosetti-Brezis-Cerami94,
Ambrosetti-GarciaAzorero-Peral00,Ambrosetti-Rabinowitz73,
Chabrowski-doO02,CortazarElguetaFelmer93,De
Figueiredo-Gossez-Ubilla03,DiazHernandez99,Willem96,Yau82}.

%
%
%

\bigskip

In \textbf{\S \ref{sec:Some Examples}}, we study particular
problems related to Einstein manifolds. Deep studies about
Einstein manifolds can be found in the books \cite{Besse87,
KramerStephaniHerltMacCallum80} and the reviews
\cite{Bourguignon84,Kazdan83,Kazdan85,Yau82}. Besides, in
\cite{Besse87} there is an approach to the existence of Einstein
warped products (see also \cite{Kim-Kim01}).

Here, we consider suitable conditions that allow us to deal with
some particular cases of the problem

\medskip
\begin{center}
\begin{minipage}{10 cm}
  \item[\textbf{(Pb-Eins.)}] Given $\mu \in \mathbb{R}$, does there exist
$\psi \in C^\infty_{>0}(B)$ such that the corresponding
\bcwpar{\psi,\mu} is an Einstein manifold?
\end{minipage}
\end{center}
\medskip

More precisely, when $B$ is an
interval in $\mathbb R$
(eventually $\mathbb{R}$) we reduce the problem to a single
ordinary differential equation that can be solved by applying
special functions. We give a more complete description if $B =
(B_m,g_B)$ is a compact scalar flat manifold, in particular when
$m=1$. Furthermore we characterize Einstein manifolds with a
precise type of metric of 2-dimensional base, generalizing
\eqref{eq:Schwarzschild to bcwp-2}. The latter result is very
close to the work of H.-J. Schmidt in \cite{Schmidt97}.

\bigskip

In the \textbf{Appendix \ref{sec:appendix}}, we give a group of
useful results about the behavior of the Laplace-Beltrami operator
under a conformal change in the metric and we present the sketch
of an alternative proof of \textbf{Theorem \ref{thm:scurv conf
warped m ge 2}} by applying a conformal change metric technique
like in \cite{DobarroLamiDozo87}.

\bigskip

\begin{center}
\textsc{Acknowledgements}
\end{center}
The authors wish to thank Diego Mazzitelli and Carmen
N\'{u}\~{n}ez for fruitful discussions of some aspects in
cosmology and string theory. F. D. was partially supported by
funds of the National Group `Analisi Reale' of the Italian
Ministry of University and Scientific Research at the University
of Trieste, he also thanks The Abdus Salam International Centre of
Theoretical Physics for their warm hospitality where part of this
work has been done.

\section{Some families of differential operators}
\label{sec:Some differential operators}

Throughout this section, $N = (N_n,h)$ is assumed to be a
pseudo-Rieman\-nian manifold of dimension $n$,
%
$|\nabla (\cdot)|^{2}=|\nabla^{N} (\cdot)|_{N}^{2}=
h(\nabla^{N}(\cdot),\nabla^{N}(\cdot))$ and
$\Delta_{h}=\Delta_{N}$.

\begin{lem}
\label{lem:2.1}
    Let $L_{h}$ be the differential operator on $C^\infty_{>0}(N)$
defined by
    \begin{equation}
     L_{h} v =\sum r_{i}\frac{\Delta_{h}v^{a_{i}}}{v^{a_{i}}},
        \label{eq:2.1}
    \end{equation}
    where any $r_{i},a_{i} \in \mathbb{R}$, $\zeta :=\sum r_{i}a_{i}
\neq 0
    $,
    $\eta:= \sum r_{i}a_{i}^{2} \neq 0 $ and the indices
    extend from $1$ to $l \in \N $. Then for $\alpha = \displaystyle
\frac{\zeta}{\eta}
    $ and
    $\beta = \displaystyle \frac{\zeta^{2}}{\eta}
    $ there
    results
    \begin{equation}
     L_{h} v= \beta
\frac{\Delta_{h}v^{\frac{1}{\alpha}}}{v^{\frac{1}{\alpha}}}.
        \label{eq:2.2}
    \end{equation}
\end{lem}

\begin{proof}
In general, for a given real value $t,$
\begin{equation}
\begin{split}
 &\nabla v^{t} = t  v^{t -1} \nabla v, \\
 &\Delta _{h} v^{t} = t [(t-1) v^{t-2} |\nabla v|^{2} + v^{t -1}
\Delta _{h}v] \textrm{ and }\\
 &\frac{\Delta _{h} v^{t}}{v^{t}} = t \left[(t-1) \frac{|\nabla
v|^{2}}{v^{2}} + \frac{\Delta
 _{h}v}{v}\right].
 \end{split}
 \label{eq:lappotencia}
\end{equation}

Thus, the right hand side of \eqref{eq:2.2}
\begin{displaymath}
\begin{split}
\beta\frac{\Delta_{h}v^{\frac{1}{\alpha}}}{v^{\frac{1}{\alpha}}}&
= \beta \frac{1}{\alpha} [(\frac{1}{\alpha} -1)\frac{|\nabla
v|^{2}}{v^{2}}+ \frac{\Delta
 _{h}v}{v}]\\
 &= \sum r_{i}a_{i}  [(\frac{\sum r_{i}a_{i}^{2}}{\sum r_{i}a_{i}}
-1)\frac{|\nabla
v|^{2}}{v^{2}}+ \frac{\Delta _{h}v}{v}] \\
&= \frac{|\nabla v|^{2}}{v^{2}} \sum_{i=1}^{l} r_{i} a_{i} (a_{i}
- 1) +
     \frac{\Delta_{h}v}{v} \sum_{i=1}^{l} r_{i} a_{i}.
 \end{split}
\end{displaymath}
And, again by \eqref{eq:lappotencia}, the left hand side of
\eqref{eq:2.2}
\begin{equation}
     L_{h} v =  \frac{|\nabla v|^{2}}{v^{2}} \sum_{i=1}^{l} r_{i}
a_{i} (a_{i} - 1) +
     \frac{\Delta_{h}v}{v} \sum_{i=1}^{l} r_{i} a_{i}.
\label{eq:operator}
\end{equation}
\end{proof}

\begin{rem}
 \label{rem:operator lemma}
 Note that equation \eqref{eq:operator} is independent of the
hypothesis $\zeta := \sum r_{i}a_{i} \neq 0 $ and
 $\eta:=\sum r_{i}a_{i}^{2} \neq 0 $, it only depends on the
 structure of the operator $L$. Thus, the following expression is
 always satisfied
 \begin{equation}
     L_{h} v =  (\eta - \zeta)\frac{|\nabla v|^{2}}{v^{2}}
     +
     \zeta \frac{\Delta_{h}v}{v}.
\label{eq:operator zetz eta}
\end{equation}
\end{rem}

\begin{cor}
\label{cor:2.1}
    Let $L_{h}$ be a differential operator defined by
    \begin{equation}
     L_{h} v = r_{1}\frac{\Delta_{h}v^{a_{1}}}{v^{a_{1}}} +
           r_{2}\frac{\Delta_{h}v^{a_{2}}}{v^{a_{2}}}  \textrm{ for }  v \in
           C_{>0}^\infty(N),
        \label{eq:2.3}
    \end{equation}
    where $ r_{1}a_{1} +
r_{2}a_{2}\neq 0 $ and
    $ r_{1}a_{1}^{2} + r_{2}a_{2}^{2}\neq 0 $. Then, by changing
the variables $v = u^{\alpha}$ with $0 < u \in
    C^\infty(N)$,
    $\alpha = \displaystyle
    \frac{ r_{1}a_{1}+r_{2}a_{2}}
    { r_{1}a_{1}^{2} + r_{2}a_{2}^{2}}$ and
    $\beta = \displaystyle
    \frac{(r_{1}a_{1} + r_{2}a_{2})^{2}}
    {r_{1}a_{1}^{2} + r_{2}a_{2}^{2}}=\alpha (r_{1}a_{1} +
r_{2}a_{2})$ there results
    \begin{equation}
     L_{h} v= \beta \frac{\Delta_{h}u}{u}.
    \label{eq:2.4}
    \end{equation}
\end{cor}

\begin{rem}
 \label{rem:Lelong-Ferrand}
To the best of our knowledge, the
only reference of an application of the identity in the form of
\eqref{eq:2.4} is an article where J. Lelong-Ferrand completed the
solution given in another paper of her about a conjecture of A.
Lichnerowicz concerning the conformal group of diffeomorphisms of
a compact $C^\infty$ Riemannian manifold, namely
if such a manifold has the group of conformal transformations,
then the manifold is globally conformal to the standard sphere of
the same dimension.
%
%
%
%
Her application corresponds to the values $r_1=1/(n-1),
r_2=-1/(n+2), a_1=n-1$ and $a_2=n$ (see \cite[p. 94 Proposition
2.2]{Lelong-Ferrand76} ).
%
%
%
%
\end{rem}

\begin{rem}
\label{rem:2.1} By the change of variables as in Corollary
\ref{cor:2.1} equations of the type
\begin{equation}
L_{h} v = r_{1}\frac{\Delta_{h}v^{a_{1}}}{v^{a_{1}}} +
           r_{2}\frac{\Delta_{h}v^{a_{2}}}{v^{a_{2}}} = H(v,x,s),
    \label{eq:nonlinear eq}
\end{equation}
transform into
\begin{equation}
\beta \Delta_{h} u = u H(u^{\alpha},x,s).
\end{equation}
We will apply this argument several times throughout the paper.
\end{rem}

\begin{exam}
\label{exam:wpsc} As it was mentioned in \S 1, the relation
connecting the scalar curvatures of the base and the fiber in a
warped product (see \cite{Beem-Ehrlich81,Beem-Ehrlich-Easley96,
Besse87,ONeil83}) is
\begin{equation}
S=-2k \frac{\Delta_{g_{B}}w}{w} - k(k-1)\frac{|\nabla^{B}
w|_{B}^{2}}{w^{2}} + S_{g_{B}} + \frac{S_{g_{F}}}{w^{2}}.
    \label{eq:scalar curv wp1}
\end{equation}
By applying \eqref{eq:lappotencia} with $t=k$ and $h=g_{B}$, it
results the following
\begin{equation}
k \frac{\Delta_{g_{B}}w}{w} + \frac{\Delta_{g_{B}}w^{k}}{w^{k}} =
-S + S_{g_{B}} + \frac{S_{g_{F}}}{w^{2}}.
    \label{eq:scalar curv wp2}
\end{equation}
Thus, by Remark \ref{rem:2.1} with $\displaystyle \alpha =
\frac{2}{k+1}$, $\displaystyle \beta = \frac{4k}{k+1}$ and $w =
u^{\alpha}$, we transform \eqref{eq:scalar curv wp1} into
\begin{equation}
\frac{4k}{k+1} \Delta_{g_{B}} u = u\left(-S + S_{g_{B}} +
\frac{S_{g_{F}}}{u^{\frac{4}{k+1}}}\right),
    \label{eq:scalar curv wp3}
\end{equation}
which is equivalent to equation \eqref{eq:scalar curv wp}
introduced in  \cite{{Dobarro87},{DobarroLamiDozo87}}.
\end{exam}

\begin{rem}
We have already mentioned that operators like $L_{h}$ are present
in different fields in \S 1. For instance, a similar situation to
\emph{Example} \ref{exam:wpsc} can be found in the study of
special cases of the Grad-Shafranov equation with a flow in plasma
physics, see \cite{
Hameiri83, Yoshida-Mahajan-Ohsaki04}.


\end{rem}


\bigskip


Now, we consider $H_{h}^{v}$ the Hessian of a function $v \in
C^\infty(N)$, so that its second covariant differential
$H_{h}^{v}=\nabla(\nabla v)$. Recall that it is the symmetric
$(0,2)$ tensor field such that for any $X,Y$ smooth vector fields
on $N$,
\begin{equation}
 \label{eq:hessian.0}
 H_{h}^{v} (X,Y) = X Y v - (\nabla_{X}Y)v = h(\nabla_{X}(\grad v),Y).
\end{equation}
Hence, for any $v \in  C^\infty_{>0}(N)$ and for all $t \in
\mathbb{R}$
\begin{equation}
 \label{eq:hessian.00}
 H_{h}^{v^{t}}  = t[(t-1)v^{t-2} dv \otimes dv  + v^{t-1}
 H_{h}^{v}],
\end{equation}
or equivalently
\begin{equation}
 \label{eq:hessian.000}
 \frac{1}{v^{t}} H_{h}^{v^{t}}  = t\left[(t-1)\frac{1}{v^{2}} dv
\otimes dv  + \frac{1}{v}
 H_{h}^{v}\right]
 ,
\end{equation}
where $\otimes$ is the usual tensorial product. Note the analogy
of the latter expressions with \eqref{eq:lappotencia} (for deeper
information about the Hessian, see p. 86 of \cite{ONeil83}).

Thus, by using the same technique applied in the proof of
\textbf{Lemma} \ref{lem:2.1} and \textbf{Remark} \ref{rem:operator
lemma}, there results
\begin{lem}
\label{lem:hessian}
 Let $\mathcal{H}_{h}$ be a differential operator on
$C^\infty_{>0}(N)$ defined by
    \begin{equation}
     \mathcal{H}_{h} v =\sum r_{i}\frac{H_{h}^{v^{a_{i}}}  }
{v^{a_{i}}},
        \label{eq:hessian.1}
    \end{equation}
    $\zeta :=\sum r_{i}a_{i}$ and
    $\eta:= \sum r_{i}a_{i}^{2}$, where the indices extend from
    $1$ to $l \in \N $ and any $r_{i},a_{i} \in \mathbb{R}$. Hence,
  \begin{equation}
  \label{eq:hessian.2}
  \mathcal{H}_{h} v  = (\eta-\zeta)\frac{1}{v^{2}} dv \otimes dv  +
\zeta \frac{1}{v}
  H_{h}^{v}.
  \end{equation}
  If furthermore, $\zeta \neq 0$ and $\eta \neq 0$, then
    \begin{equation}
     \mathcal{H}_{h} v= \beta \frac{H_{h}^{v^{\frac{1}{\alpha}}}
}
     {v^{\frac{1}{\alpha}}},
        \label{eq:hessian.3}
    \end{equation}
    where $\alpha = \displaystyle \frac{\zeta}{\eta}$ and
    $\beta = \displaystyle \frac{\zeta^{2}}{\eta}$.

\end{lem}

%

\section{About base conformal warped products}
\label{sec:base conformal warped products}


In this section, we define precisely {\it base conformal warped
products} and compute covariant derivatives and curvatures of {\it
base conformal warped products}. Several proofs contain standard
but long computations, and hence will be omitted.

Let $(B,g_B)$ and  $(F,g_F)$ be $m$ and $k$ dimensional {\it
pseudo-Riemannian manifolds}, respectively. Then  $M=B \times F$
is an $(m+k)$-dimensional  {\it pseudo-Riemannian manifold} with
$\pi\colon B \times F \to B$ and $\sigma \colon B \times F \to F$
the usual projection maps.

Throughout this paper we use the {\it natural product coordinate
system} on the {\it product manifold} $B \times F$, namely. Let
$(p_0,q_0)$ be a point in $M$ and {\it coordinate charts} $(U,x)$
and $(V,y)$ on $B$ and $F$, respectively such that $p_0 \in B$ and
$q_0 \in F$. Then we can define a {\it coordinate chart} $(W,z)$
on $M$ such that $W$ is an open subset in $M$ contained in $U
\times V$, $(p_0,q_0) \in W$ and for all $(p,q)$ in $W$,
$z(p,q)=(x(p),y(q)),$ where $x=(x^1,\cdots, x^m)$ and
$y=(y^{m+1},\cdots,y^{m+k}).$

Clearly, the set of all $(W,z)$ defines an atlas on $B \times F$.
Here, for our convenience, we call the $j$-th {\it component} of
$y$ as $y^{m+j}$ for all $j \in \{1, \cdots, k \}$.

Let $\phi \colon B\to \mathbb R \in {\mathcal C}^\infty(B)$ then
the lift of $\phi$ to $B \times F$ is $\widetilde{\phi}= \phi
\circ \pi \in {\mathcal C}^\infty(B \times F),$ where ${\mathcal
C}^\infty(B)$ is the set of all smooth real-valued functions on
$B.$

Moreover, one can define {\it lifts} of {\it tangent vectors} as:
Let $X_p \in T_p(B)$ and $q \in F$ then the lift
$\widetilde{X}_{(p,q)}$ of $X_p$ is the unique tangent vector in $
T_{(p,q)}(B \times \{q \})$ such that $d
\pi_{(p,q)}(\widetilde{X}_{(p,q)}) =X_p$ and  $d \sigma_{(p,q)}
(\widetilde{X}_{(p,q)})=0.$ We will denote the set of all lifts of
all tangent vectors of $B$ by $L_{(p,q)}(B).$

Similarly, we can define {\it lifts} of {\it vector fields}. Let
$X \in \mathfrak X(B)$ then the lift of $X$ to $B \times F$ is the
vector field $\widetilde X \in \mathfrak {X}(B \times F)$ whose
value at each $(p,q)$ is the lift of $X_p$ to $(p,q).$ We will
denote the set of all lifts of all vector fields of $B$ by
$\mathfrak L(B).$

\begin{dfn} Let $(B,g_B)$ and $(F,g_F)$ be {\it pseudo-Riemannian}
manifolds and also let $w \colon B \to (0,\infty)$ and $c \colon B
\to (0,\infty)$ be smooth functions. The {\it base conformal
warped product} (briefly \textit{bcwp}) is the {\it product
manifold} $B \times F$ furnished with the metric tensor
$g=c^{2}g_B \oplus w^{2}g_F$ defined by
\begin{equation} \label{dwp}
g=(c \circ \pi)^2 \pi^{\ast}(g_B) \oplus (w \circ \pi)^2
\sigma^{\ast}(g_F).
\end{equation}
By analogy with \cite{PongeReckziegel93} we will denote this
structure by \twistpar{B}{F}{c}{w}. The function $w \colon B \to
(0,\infty)$ is called the \textit{warping function} and the function
$c \colon B \to (0,\infty)$ is said to be the \textit{conformal
factor}.

If $c \equiv 1$ and $w$ is not identically $1$, then we obtain a
{\it singly warped product}. If both $w \equiv 1$ and $c \equiv
1,$ then we have a {\it product manifold}. If neither $w$ nor $c$
is constant, then we have a {\it nontrivial bcwp}.

If $(B,g_B)$ and $(F,g_F)$ are both {\it Riemannian manifolds,}
then \twistpar{B}{F}{c}{w}
is also a {\it Riemannian manifold}. We call \twistpar{B}{F}{c}{w}
as a
{\it Lorentzian base
conformal warped product} if $(F,g_F)$ is {\it Riemannian} and
either $(B,g_B)$ is {\it Lorentzian} or else $(B,g_B)$ is a
one-dimensional manifold with a {\it negative definite} metric
$-dt^2$.
\end{dfn}

\begin{nota}\label{nota: bcwp1}
From now on, we will identify the operators defined on the base
(respectively, fiber) of a \bcwpp with the name of the base
(respectively, fiber) as a sub or super index. Unlike, the
operators defined on the whole \bcwpp
will not have labels. For instance, the Riemann curvature tensor
of the base $(B,g_B)$ will be denoted by $R_B$ and likewise $R_F$
denotes for that of the fiber $(F,g_F).$ Thus, the Riemann
curvature tensor of \twistpar{B}{F}{c}{w} is denoted by $R$.
\end{nota}

\subsection{Covariant Derivatives}
\label{subsec:Covariant Derivatives}

We state the covariant derivative formulas and the geodesic
equation for a base conformal warped product manifold
\twistpar{B}{F}{c}{w}.

The gradient operator of smooth functions on \twistpar{B}{F}{c}{w}
is denoted by $\nabla$ and $\nabla^B$ and $\nabla^F$ denote the
gradients of $(B,g_B)$ and $(F,g_F)$, respectively (see
\textit{Notation \ref{nota: bcwp1}}).

\begin{prp} \label{grad1} Let $\phi \in {\mathcal C}^\infty(B)$
and $\psi \in {\mathcal C}^\infty(F).$ Then
$$ \nabla\phi = \frac{1}{c^2}\nabla^B \phi \quad \text{and}
\quad \nabla\psi = \frac{1}{w^2}\nabla^F \psi.$$
\end{prp}

Also, we express the {\it covariant derivative} on $B \times F$ in
terms of the {\it covariant derivatives} on $B$ and $F$ by using
the Kozsul formula, which takes the following form on a base
conformal warped product as above: Let $X,Y,Z \in \mathfrak L(B)$
and $V,W,U \in \mathfrak L(F),$ then

\begin{eqnarray*} 2g(\nabla_{X+V}(Y+W),Z+U)
& = & (X+V)g(Y+W,Z+U) \\
& + & (Y+W)g(X+V,Z+U) \\
& - & (Z+U)g(X+V,Y+W) \\
& + & g([X+V,Y+W],Z+U) \\
& - & g([X+V,Z+U],Y+W) \\
& - & g([Y+W,Z+U],X+V),
\end{eqnarray*}
where $[\cdot,\cdot]$ denotes the Lie bracket.
\begin{thm} \label{cov} Let $X,Y \in \mathfrak L(B)$
and $V,W \in \mathfrak L(F)$. Then
\begin{enumerate}
\item ${\displaystyle \nabla_X Y = \nabla_X^B Y + \frac{X( c)}{
c}Y+ \frac{Y( c)}{ c}X-\frac{g_B(X,Y)}{ c}\nabla^B c}$,
 \item
${\displaystyle \nabla_X V = \nabla_V X=\frac{X( w)}{ w}V}$,
 \item
${\displaystyle \nabla_V W = \nabla_V^F W-\frac{ w}{ c^2} g_F(V,W)
\nabla^B w}$.
\end{enumerate}
\end{thm}

\begin{rem} Let $X,Y \in \mathfrak L(B)$
and $V,W \in \mathfrak L(F)$. If $[\cdot,\cdot]$ denotes for the
Lie bracket on \twistpar{B}{F}{c}{w},
then $\displaystyle{[X,Y]=[X,Y]_B},$ \, $\displaystyle{[X,V]=0}$
and $\displaystyle{[V,W]=}$ $[V,W]_F$.
\end{rem}

\begin{prp} Let $(p,q) \in$ \twistpar{B}{F}{c}{w}. Then
\begin{enumerate}
\item The leaf $B \times \{ q \}$ and the fiber $\{ p \} \times F$ are totally umbilic.
\item The leaf $B \times \{ q \}$ is
totally geodesic.
\item The fiber $\{ p \} \times F$ is totally
geodesic when $(\nabla^B w)(p)=0.$
\end{enumerate}
\end{prp}

Now, we will establish the {\it geodesic equations} for {\it base
conformal warped products}. The version for {\it singly warped
products} is well known (compare page 207 of \cite{{ONeil83}}).

\begin{prp} \label{geoeqn} Let $\gamma=(\alpha, \beta)
\colon I \to$\twistpar{B}{F}{c}{w} be a (smooth) curve where $I
\subseteq \mathbb R$. Then $\gamma=(\alpha, \beta)$ is a geodesic
in \twistpar{B}{F}{c}{w}
if and only if for any $t \in I$,
\begin{enumerate}
\item $\displaystyle{\alpha^{\prime
\prime}=-2\frac{\alpha^\prime(c)}
{c}\alpha^\prime+\frac{g_B(\alpha^\prime,\alpha^\prime)}{c}
\nabla^B c +\frac{w g_F(\beta^\prime, \beta^\prime)} {c^2}\nabla^B
w}$,
\item $\displaystyle{\beta^{\prime
\prime}=-2\frac{\alpha^\prime(w)}{w} \beta^\prime}$.
\end{enumerate}
\end{prp}

\begin{rem} If $\gamma=(\alpha, \beta) \colon I \to$\twistpar{B}{F}{c}{w}
is a geodesic in \twistpar{B}{F}{c}{w},
then $\beta \colon
I \to F$ is a pre-geodesic in $(F,g_F).$
\end{rem}

\subsection{Riemannian Curvatures}
\label{subsec:Riemannian Curvatures}

From now on, we use the definition and the sign convention for the
{\it curvature} as in \cite[p. 16-25]{Beem-Ehrlich-Easley96} (note
the difference with \cite{ONeil83}), namely. For an arbitrary
$n$-dimensional pseudo-Riemannian manifold $(N,h)$, letting $X,Y,Z
\in \mathfrak L(N)$, we take the Riemann curvature tensor
\begin{equation*}\label{}
  R(X,Y)Z=\nabla_X \nabla_Y Z - \nabla_Y \nabla_X Z -
  \nabla_{[X,Y]} Z.
\end{equation*}
Furthermore, for each $p \in N $, the Ricci curvature tensor is
given by
\begin{equation*}\label{}
  \ric (X,Y)= \sum_{i=1} ^n h(E_i,E_i) h(R(E_i,Y)X,E_i),
\end{equation*}
where $\{E_1, \cdots , E_n\}$ is an orthonormal basis for $T_p N$.

Now, we give the Riemannian curvature formulas for a base
conformal warped product. But first we state the Hessian tensor
denoted by ${\rm H}$ (see \S \ref{sec:Some differential
operators}) on this class of warped products.

\begin{prp}\label{prp: hessian bcwp}Let $X,Y \in \mathfrak L(B)$ and
$V,W \in \mathfrak L(F)$ and also let $\phi \in {\mathcal
C}^\infty(B)$ and $\psi \in {\mathcal C}^\infty(F).$ Then, the
Hessian $H$ of \twistpar{B}{F}{c}{w} satisfies
\begin{enumerate}
\item $\displaystyle{{\rm H}^\phi(X,Y)={\rm H}_B^\phi(X,Y)+
\frac{g_B(X,Y)}{c}g_B(\nabla^B \phi ,\nabla^B c )}$ \\
\mbox

\hspace{1.58 cm} $\displaystyle{-\frac{X(c)Y(\phi)}{c}-
\frac{Y(c)X(\phi)}{c}},$
\item $\displaystyle{{\rm
H}^\psi(X,Y)=0},$
\item $\displaystyle{{\rm H}^\phi(X,V)=0},$
\item
$\displaystyle{{\rm H}^\psi(X,V)=-\frac{X(w)V(\psi)}{w}},$
\item
$\displaystyle{{\rm H}^\phi(V,W)=\frac{w}{c^2}
g_F(V,W)g_B(\nabla^B w ,\nabla^B \phi )},$
\item
$\displaystyle{{\rm H}^\psi(V,W)={\rm H}_F^\psi(V,W)}.$
\end{enumerate}
\end{prp}

\begin{thm} \label{riemc} Let $X,Y,Z \in \mathfrak L(B)$ and
$V,W,U  \in \mathfrak L(F).$ Then, the curvature Riemann tensor
$R$ of \twistpar{B}{F}{c}{w} satisfies
\begin{enumerate}
\item $\displaystyle{{\rm R}(X,Y)Z={\rm R}_B(X,Y)Z-\frac
{{\rm H}^c(Y,Z)}{c}X+\frac{{\rm H}^c(X,Z)}{c}Y}$ \\
\mbox

\hspace{1.68 cm} $\displaystyle{+2\frac{X(c)}{c^2}g_B(Y,Z)
\nabla^B c -2\frac{Y(c)}{c^2}g_B(X,Z) \nabla^B c }$ \\
\mbox

\hspace{1.68 cm} $\displaystyle{+\frac{g_B(X,Z)}{c}\nabla^B_Y
\nabla^B c -\frac{g_B(Y,Z)}{c}\nabla^B_X \nabla^B c },$
\item $\displaystyle{{\rm R}(X,V)Y=\frac{{\rm H}^w(X,Y)}{w}V},$
\item $\displaystyle{{\rm R}(X,Y)V=0},$
\item $\displaystyle{{\rm R}(V,W)X=0},$
\item $\displaystyle{{\rm R}(V,X)W=w g_F(V,W){\rm h}^w(X)}$,
\item $\displaystyle{{\rm R}(V,W)U={\rm R}_F(V,W)U}$ \\
\mbox

\hspace{1.68 cm} $\displaystyle{+\frac {g_B(\nabla^B w , \nabla^B
w )}{c^2} \Bigl(g_F(V,U)W-g_F(W,U)V \Bigl)}$,
\end{enumerate}
where ${\rm h}^w(X)$ is given in the remark that follows.
\end{thm}

\begin{rem} \label{rem: hw} Note that $\displaystyle{{\rm h}^w(X)=
\nabla_X \nabla w }$ and $\displaystyle{\nabla w = \frac{1}{c^2}
\nabla^B w }.$ Hence,
\begin{eqnarray*}
{\rm h}^w(X) & = & -2\frac{X(c)}{c^3} \nabla^B w +\frac{1}{c^2}
\left(\nabla_X^B \nabla^B w +\frac{X(c)}{c}\nabla^B w  \right.\\
& + & \left. \frac{g_B(\nabla^B w ,\nabla^B c )}{c}X -
\frac{X(w)}{c} \nabla^B c  \right).
\end{eqnarray*}
\end{rem}

\subsection{Ricci Curvatures}
\label{subsec:Ricci Curvatures}

We compute Ricci curvatures of the base conformal warped product
applying that if $\{E_1, ... , E_m\}$ is a $g_B-$orthonormal frame
field on an open set $U \subseteq B$ and $\{\tilde{E}_{m+1} ....
\tilde{E}_{m+k}\}$ is a $g_F-$orthonormal frame field on an open
set $V \subseteq F$, then
$$
\{c^{-1}E_1, ... , c^{-1}E_m,w^{-1}\tilde{E}_{m+1} ....
w^{-1}\tilde{E}_{m+k}\}
$$
is a $g-$orthonormal frame field on an open set $W \subseteq U
\times V \subseteq B \times F$.

\begin{prp} \label{prp: laplacian bcwp}Let $\phi \in {\mathcal
C}^\infty(B)$ and $\psi \in {\mathcal C}^\infty(F).$ Then, the
Laplace-Beltrami operator $\Delta$ of \twistpar{B}{F}{c}{w}
satisfies
\begin{enumerate}
\item $\displaystyle{\Delta \phi =\frac{\Delta_B \phi }{c^2}+
\frac{m-2}{c^3}g_B(\nabla^B \phi ,\nabla^B c )} +
\frac{1}{c^2} \frac{k}{w} g_B(\nabla^B w ,\nabla^B \phi )$,
\item $\displaystyle{\Delta \psi =\frac{\Delta_F \psi }{w^2}}$.
\end{enumerate}
\end{prp}

\begin{thm} \label{riki} Let $X,Y \in \mathfrak L(B)$ and
$V,W \in \mathfrak L(F)$. Then, the Ricci tensor $Ric$ of
\twistpar{B}{F}{c}{w} satisfies


\begin{enumerate}
\item $\displaystyle{{\rm Ric}(X,Y)={\rm Ric}_B(X,Y)}$\\
\mbox

\hspace{1.67 cm} $\displaystyle{ - (m-2)\frac{1}{c}
{\rm H}_{B}^c(X,Y)+2(m-2)\frac{1}{c^{2}}X(c)Y(c)}$\\
\mbox

\hspace{1.67 cm} $\displaystyle -\left[(m-3)\frac{g_B(\nabla^B c
 ,\nabla^B c )}{c^{2}}+\frac{\Delta_{B}c}{c}
\right]g_B(X,Y)$\\
\mbox

\hspace{1.67 cm} $\displaystyle{-k\frac{1}{w} {\rm H}_{B}^w(X,Y) -
k
\frac{g_B(\nabla^B w ,\nabla^B c )}{wc}g_B(X,Y)}$ \\
\mbox

\hspace{1.67 cm} $\displaystyle
+k\frac{X(c)}{c}\frac{Y(w)}{w}+k\frac{Y(c)}{c}\frac{X(w)}{w}$,
\mbox

\item $\displaystyle{{\rm Ric}(X,V)=0}$,

\item $\displaystyle{{\rm Ric}(V,W)={\rm Ric}_F(V,W)}$ \\
\mbox

\hspace{1.67 cm} $\displaystyle -\frac{w^{2}}{c^{2}}g_F(V,W)
\left[ (m-2)
\frac{g_B(\nabla^B w ,\nabla^B c )}{wc}
+\frac{\Delta_{B}w}{w}\right.$\\
\mbox

\hspace{1.67 cm} $\displaystyle{\left. + (k-1) \frac{g_B(\nabla^B
w ,\nabla^B w )}{w^2}\right]}$.
\end{enumerate}

\end{thm}

An equivalent formulation of \textbf{Theorem \ref{riki}} is

\begin{thm} \label{global riki}
The Ricci tensor $Ric$ of \twistpar{B}{F}{c}{w} satisfies

\begin{enumerate}
\item $\displaystyle{{\rm Ric}={\rm Ric}_B-
\left[(m-2)\frac{1}{c}{\rm H}_{B}^c + k\frac{1}{w} {\rm H}_{B}^w
\right]}$\\
\mbox

\hspace{0.67 cm}  $\displaystyle{ +2(m-2)\frac{1}{c^{2}} dc
\otimes dc +k\frac{1}{wc} [dc \otimes dw + dw
\otimes dc]}$\\
\mbox

\hspace{0.67 cm} $\displaystyle -\left[(m-3)\frac{g_B(\nabla^B c
 ,\nabla^B c )}{c^{2}}+\frac{\Delta_{B}c}{c} + k
\frac{g_B(\nabla^B w ,\nabla^B c )}{wc}
\right]g_B $\\
\mbox

\hspace{0.67 cm} $\displaystyle \textrm{ on } \mathcal{L}(B)\times
\mathcal{L}(B)$, \mbox

\item $\displaystyle{{\rm Ric}=0} \textrm{ on }
\mathcal{L}(B)\times \mathcal{L}(F)$,

\item $\displaystyle{{\rm Ric}={\rm Ric}_F
-\frac{w^{2}}{c^{2}} \left[ (m-2)
\frac{g_B(\nabla^B w ,\nabla^B c )}{wc} +
\frac{\Delta_{B}w}{w}\right.}$ \\
\mbox

\hspace{0.67 cm} $\displaystyle{\left.+ (k-1) \frac{g_B(\nabla^B w
,\nabla^B w )}{w^2}\right]  g_F} \textrm{ on }
\mathcal{L}(F)\times \mathcal{L}(F)$.
\end{enumerate}
\end{thm}

\begin{rem} If $m \neq 2$ and $ k \neq 1$, applying
\eqref{eq:lappotencia}, the expression of the Ricci tensor of
\twistpar{B}{F}{c}{w} in \textbf{Theorem \ref{riki}} may be
written as
\begin{enumerate}
\item $\displaystyle{{\rm Ric}(X,Y)={\rm Ric}_B(X,Y)}$\\
\mbox

\hspace{1.67 cm} $\displaystyle{ - (m-2)\frac{1}{c}
{\rm H}_{B}^c(X,Y)+2(m-2)\frac{1}{c^{2}}X(c)Y(c)}$\\
\mbox

\hspace{1.67 cm} $\displaystyle
-\frac{1}{m-2}\frac{\Delta_{B}c^{m-2}}{c^{m-2}}
g_B(X,Y)$\\
\mbox

\hspace{1.67 cm} $\displaystyle{-k\frac{1}{w} {\rm H}_{B}^w(X,Y) -
k \frac{1}{wc} g_B(
\nabla^B w ,\nabla^B c )g_B(X,Y)}$ \\
\mbox

\hspace{1.67 cm} $\displaystyle
+k\frac{X(c)}{c}\frac{Y(w)}{w}+k\frac{Y(c)}{c}\frac{X(w)}{w}$,
\mbox

\item $\displaystyle{{\rm Ric}(X,V)=0}$,

\item $\displaystyle{{\rm Ric}(V,W)={\rm Ric}_F(V,W)}$ \\
\mbox

\hspace{1.67 cm} $\displaystyle-\frac{w^{2}}{c^{2}}g_F(V,W) \left[
(m-2)\frac{1}{wc}  g_B(\nabla^B w ,\nabla^B c ) +
\frac{1}{k}\frac{\Delta_{B}w^{k}}{w^{k}}\right]$.
\end{enumerate}
\end{rem}

\bigskip

\subsection{Scalar Curvature}
\label{subsec:Scalar Curvature}

By using the orthonormal frame introduced above, one can obtain
the following result after a standard computation.

\begin{thm}\label{sca-c} The scalar curvature $S$ of
\twistpar{B}{F}{c}{w} is given by
\begin{eqnarray*}
c^{2}S & = & S_B + S_F\frac{c^{2}}{w^2}
-2(m-1)\frac{\Delta_{B}c}{c}- 2k
\frac{\Delta_{B} w}{w}\\
& - & (m-4)(m-1)   \frac{g_B(\nabla^B c ,\nabla^B c )}{c^{2}}\\
& - & 2k(m-2) \frac{g_B(\nabla^B w ,\nabla^B c )}{w c} \\
& - & k(k-1)\frac{g_B(\nabla^B w ,\nabla^B w )}{w^{2}}.
\end{eqnarray*}
\end{thm}

\section{Curvature of $({(B \times F)}_{m+k},\psi^{2\mu} g_{B} +
\psi^{2}g_{F})$} \label{sec:Curvature of sbcwp}

From now on, we will deal with \bcwp's, i.e.
\twistpar{B}{F}{\psi^\mu}{\psi}, and specifically concentrate on
its Ricci tensor and scalar curvature.

\subsection{Ricci Tensor}\label{subsection sbcwp Ricci}

\begin{thm}
\label{thm:ricci sbcwp m ge 3} Let $B = (B_m,g_{B})$ and $F =
(F_k,g_{F})$ be two pseudo-Rieman\-nian manifolds with dimensions
$m \geq 3$ and $k \geq 1$, $\displaystyle \mu \in
\mathbb{R}\setminus \{0, 1, \overline{\mu},
\overline{\mu}_{\pm}\}$ with
$\overline{\mu}:=\displaystyle-\frac{k}{m-2}$ and
$\overline{\mu}_{\pm}:=\overline{\mu}\pm\sqrt{\overline{\mu}^{2}-\overline{\mu}}$
and $\psi \in C^\infty_{>0}(B)$. Then, the Ricci curvature tensor
\ric ~ of the base conformal warped product
\twistpar{B}{F}{\psi^\mu}{\psi}
verifies the relation

\begin{equation}\label{eq:ricci special 1}
\begin{split}
&\displaystyle{{\rm Ric}={\rm Ric}_B+
\beta^{H}\frac{1}{\psi^{\frac{1}{\alpha ^{H}}}}{\rm
H}_{B}^{\psi^{\frac{1}{\alpha ^{H}}}} - \beta^{\Delta}
\frac{1}{\psi ^{\frac{1}{\alpha ^{\Delta}}}} {\Delta}_{B}
\psi^{\frac{1}{\alpha ^{\Delta}}} }g_B
\textrm{ on } \mathcal{L}(B)\times \mathcal{L}(B),
\\
&\displaystyle{{\rm Ric}=0} \textrm{ on } \mathcal{L}(B)\times
\mathcal{L}(F),
\\
&\displaystyle{{\rm Ric}={\rm Ric}_F
-\frac{1}{ \psi^{2(\mu - 1)}}
\frac{\beta^{\Delta}}{\mu} \frac{1}{\psi ^{\frac{1}{\alpha
^{\Delta}}}} {\Delta}_{B} \psi^{\frac{1}{\alpha ^{\Delta}}}
g_F}
\textrm{ on } \mathcal{L}(F)\times \mathcal{L}(F),
\end{split}
\end{equation}
where
\begin{equation}\label{eq:alphaDH betaDH particular}
  \begin{array}{lcl}
    \alpha^{\Delta} &=&   \displaystyle \frac{1}{(m-2)\mu + k},      \\
    \beta^{\Delta} &=& \displaystyle \frac{\mu}{(m-2)\mu + k},      \\
    \alpha^{H} &=& \displaystyle \frac{-[(m-2)\mu + k]}{\mu [(m-2)\mu
+ k]+k(\mu - 1)},    \\
    \beta^{H} &=&  \displaystyle \frac{[(m-2)\mu + k]^{2}}{\mu
[(m-2)\mu + k]+k(\mu - 1)}.
  \end{array}
\end{equation}

\end{thm}

\begin{proof} Applying \textbf{Theorem \ref{global riki}} with
$c=\psi^{\mu}$ and $w=\psi$, we obtain
\begin{equation}\label{eq:ricci special 2}
\begin{split}
\displaystyle \ric = & \ric_B -
\left[(m-2)\frac{1}{\psi^{\mu}}{\rm H}_{B}^{\psi^{\mu}} +
k\frac{1}{\psi} {\rm H}_{B}^\psi \right]
\\
&\displaystyle +2\mu[(m-2)\mu +k ] \frac{1}{\psi^{2}}d\psi \otimes
d\psi
\\
&\displaystyle -\left[((m-3)\mu^{2}+k\mu)\frac{g_B(\nabla^B \psi
,\nabla^B \psi
)}{\psi^{2}}+\frac{\Delta_{B}\psi^{\mu}}{\psi^{\mu}} \right]g_B
\\
&\displaystyle \textrm{ on } \mathcal{L}(B)\times \mathcal{L}(B),
\\
\ric = & 0 \textrm{ on } \mathcal{L}(B)\times \mathcal{L}(F),
\\
\displaystyle \ric =& \ric_F
-\frac{1}{\psi^{2(\mu - 1)}} \left[ ((m-2)\mu + k - 1)
\frac{g_B(\nabla^B \psi,\nabla^B \psi)}{\psi^{2}} +
\frac{\Delta_{B} \psi}{\psi} \right]  g_F
\\
&\textrm{ on } \mathcal{L}(F)\times \mathcal{L}(F).
\\
\end{split}
\end{equation}
So by \eqref{eq:hessian.000} and  \eqref{eq:lappotencia}, with
$t=\mu\neq 0, 1$, there results
\begin{equation}\label{eq:ricci special 3}
\begin{split}
\displaystyle{\rm Ric}=&{\rm Ric}_B+
\left[r^{H}_{1}\frac{1}{\psi^{\mu}}{\rm H}_{B}^{\psi^{\mu}} +
r^{H}_{2}\frac{1}{\psi} {\rm H}_{B}^\psi \right]
\\
&\displaystyle
-\left[r^{\Delta}_{1}\frac{\Delta_{B}\psi^{\mu}}{\psi^{\mu}} +
r^{\Delta}_{2}\frac{\Delta_{B}\psi}{\psi} \right]g_B
\textrm{ on } \mathcal{L}(B)\times \mathcal{L}(B),
\\
\displaystyle{\rm Ric}= & 0 \textrm{ on } \mathcal{L}(B)\times
\mathcal{L}(F),
\\
\displaystyle{\rm Ric}= & {\rm Ric}_F
-\frac{1}{\psi^{2(\mu - 1)}} \left[ ((m-2)\mu + k - 1)
\frac{g_B(\nabla^B \psi,\nabla^B \psi )}{\psi^{2}} +
\frac{\Delta_{B}\psi}{\psi}\right]  g_F
\\
&\textrm{ on } \mathcal{L}(F)\times \mathcal{L}(F),
\end{split}
\end{equation}
where
\begin{equation*}
  \begin{array}{lcl}
    (\mu-1)r^{H}_{1} &=& (m-2)\mu + m-2+2k,           \\
    (\mu-1)r^{H}_{2} &=&  -(m-2)2\mu^{2} - k(3\mu -1),           \\
    (\mu-1)r^{\Delta}_{1} &=&  (m-2)\mu + k -1,     \\
    (\mu-1)r^{\Delta}_{2} &=& -\mu ((m-2)\mu + k - \mu).
  \end{array}
\end{equation*}
Hence, using the notation introduced in \textbf{Lemmas
\ref{lem:hessian}} and \textbf{\ref{lem:2.1}} and \textit{Remark}
\textit{\ref{rem:operator lemma}},
\begin{equation} \label{eq:ricci special 4}
\begin{split}
\displaystyle{\rm Ric}=&{\rm Ric}_B+
(\eta^{H} - \zeta^{H} )\frac{1}{\psi^{2}}d\psi \otimes
d\psi+\zeta^{H} \frac{1}{\psi} {\rm H}_{B}^\psi
\\
&\displaystyle -\left[ (\eta^{\Delta} - \zeta^{\Delta} )
\frac{g_B(\nabla^B \psi ,\nabla^B \psi)}{\psi^{2}}
+
\zeta^{\Delta}\frac{\Delta_{B}\psi}{\psi} \right]g_B
\textrm{ on } \mathcal{L}(B)\times \mathcal{L}(B),
\\
\displaystyle{\rm Ric}=&0 \textrm{ on } \mathcal{L}(B)\times
\mathcal{L}(F),
\\
\displaystyle{\rm Ric}=&{\rm Ric}_F
-\frac{1}{\psi^{2(\mu - 1)}} \left[
\left(\frac{\eta^{\Delta}}{\mu} - \frac{\zeta^{\Delta}}{\mu}
\right ) \frac{g_B(\nabla^B \psi ,\nabla^B \psi )}{\psi^{2}} +
\frac{\zeta^{\Delta}}{\mu}\frac{\Delta_{B}\psi}{\psi}\right] g_F
\\
&\textrm{ on } \mathcal{L}(F)\times \mathcal{L}(F),
\end{split}
\end{equation}
where
\begin{equation}\label{eq:zeta eta HD}
  \begin{array}{lclcl}
    \zeta^{H} &=& r^{H}_{1} \mu+ r^{H}_{2}&=&-[(m-2)\mu +k ],
\\
    \eta^{H} &=& r^{H}_{1} \mu^{2}+ r^{H}_{2} &=&\mu[(m-2)\mu + k] +
k(\mu-1),      \\
    \zeta^{\Delta} &=& r^{\Delta}_{1} \mu+ r^{\Delta}_{2} &=&\mu,
\\
    \eta^{\Delta} &=& r^{\Delta}_{1} \mu^{2}+ r^{\Delta}_{2}&=&\mu
[(m-2)\mu +
    k].
  \end{array}
\end{equation}
Note that
\begin{equation}\label{eq:zeta eta HD exceptional }
  \begin{array}{lcccccl}
    \zeta^{H} &=& 0 & \Leftrightarrow
&\mu&=&\overline{\mu}:=\displaystyle-\frac{k}{m-2},       \\
    \eta^{H} &=& 0 & \Leftrightarrow
&\mu&=&\overline{\mu}_{\pm}:=\overline{\mu}\pm\sqrt{\overline{\mu}^{2}-\overline{\mu}},
       \\
    \zeta^{\Delta} &=& 0 & \Leftrightarrow  &\mu &=&  0,    \\
    \eta^{\Delta} &=& 0 & \Leftrightarrow &\mu
&=&0,\displaystyle-\frac{k}{m-2}.
  \end{array}
\end{equation}
So, if $\mu \in \mathbb{R}\setminus \{0, 1, \overline{\mu},
\overline{\mu}_{\pm}
\}$ and considering
\begin{equation}\label{eq:zeta eta HD formulas}
  \begin{array}{lcl}
    \alpha^{\Delta} &=&   \displaystyle
\frac{\zeta^{\Delta}}{\eta^{\Delta}},      \\
    \beta^{\Delta} &=& \displaystyle
\frac{(\zeta^{\Delta})^{2}}{\eta^{\Delta}},  \\
    \alpha^{H} &=& \displaystyle \frac{\zeta^{H}}{\eta^{H}},    \\
    \beta^{H} &=&  \displaystyle \frac{(\zeta^{H})^{2}}{\eta^{H}},
  \end{array}
\end{equation}
along with \textbf{Lemmas \ref{lem:hessian}} and
\textbf{\ref{lem:2.1}} results the thesis.
\end{proof}

\begin{rem} We will make some comments about the previous results
and compare the above formulas with Ricci tensor formulas in the
case of a conformal manifold and a warped product.

\noindent
\begin{description}
  \item[i] Note that the system \eqref{eq:ricci special 2}
remains valid without conditions on $\mu$, $m \ge 1$ and $k \ge
0$.
  \item[ii] The system \eqref{eq:ricci special 2} with
$\mu = 1$, $m \ge 1$ and $k = 0$ give the expression of the Ricci
tensor under a conformal change in the base given by
$\tilde{g}_{B}=\psi^{2} g_{B}$, where $\psi \in C^\infty_{>0}(B)$
(see \cite{Aubin98}, \cite{Besse87}).
  \item[iii] For $\mu = 0$, $m \ge 1$ and $k \ge
1$ the system \eqref{eq:ricci special 2} reproduces the
expressions of the Ricci tensor for a singly warped product (
\cite{Beem-Ehrlich-Easley96}, \cite{Besse87}, \cite{ONeil83}).
\end{description}
\end{rem}

\bigskip

The \textsc{Table \ref{tab:ricci mu exceptional cases m ge 3}} is a
synthesis of the $\mu-$exceptional cases in the \textbf{Theorem
\ref{thm:ricci sbcwp m ge 3}}. In that table $\zeta^{H}$,
$\eta^{H}$, $\zeta^{\Delta}$ and $\eta^{\Delta}$ are computed with
the final expressions of \eqref{eq:zeta eta HD}. This is the reason
to include the column titled ``formal system", and hence the systems
written in that column are justified a posteriori.


\begin{small}
\begin{table}
\caption{\small Einstein equations, ~$m \ge 3$, ~$\mu-$exceptional
cases in \textbf{Theorem \ref{thm:ricci sbcwp m ge 3}}}
\label{tab:ricci mu exceptional cases m ge 3}
\begin{tabular}{ccccccccc}
 \hline\hline
   $\mu$ & $m$ & $k$ & $\zeta^{H}$ & $\eta^{H}$ & $\zeta^{\Delta}$ &
$\eta^{\Delta}$
   & \begin{minipage}{.9 cm}
    \textit{ \centerline{genuine} \\ \centerline{system}}
    \end{minipage}
   &  \begin{minipage}{.9 cm}
    \textit{ \centerline{formal} \\ \centerline{system}}
    \end{minipage}
\\\\
 \hline
\\
  0 & $\ge 3$ & $\ge 1 $ & $-k$ & $-k$ & $0$ & $0$ &
  \eqref{eq:ricci special 2} & \eqref{eq:ricci special 4}
  \\\\
  1 & $\ge 3$ & $\ge 1 $ & $-[m-2+k]$ & $m-2+k$ & $1$ & $m-2+k$ &
\eqref{eq:ricci special 2}& \eqref{eq:ricci special 4}
\\\\
  $\overline{\mu}$ & $\ge 3$ & $\ge 1 $ & 0 & $k(\overline{\mu}-1)$ &
$\overline{\mu}$ & 0 & \eqref{eq:ricci special 4} & -
\\\\
  $\overline{\mu}_{\pm}$& $\ge 3$ & $\ge 1 $ &
  $\displaystyle k
\frac{\overline{\mu}_{\pm}-1}{\overline{\mu}_{\pm}}$
 & 0 &
$\overline{\mu}_{\pm}$ &
$-k(\overline{\mu}_{\pm}-1)$
 &\eqref{eq:ricci special 4} & -
\\\\
  \hline\hline
\end{tabular}
\end{table}
\end{small}


\begin{rem} \label{thm:ricci sbcwp m = 1, 2} Here, we
consider the cases $m=1$ and $m=2$, with $ k\ge 1$. The results
and the proof are essentially the same as \textbf{Theorem
\ref{thm:ricci sbcwp m ge 3}}, but the conditions \eqref{eq:zeta
eta HD exceptional } take the following form.
\begin{description}
  \item[$\underline{m=1}$]
  \begin{equation}\label{eq:zeta eta HD exceptional m=1}
  \begin{array}{lcccccl}
    \zeta^{H} &=& 0 & \Leftrightarrow &\mu&=&k ,         \\
    \eta^{H} &=& 0 & \Leftrightarrow
&\mu&=&\overline{\mu}_{\pm}:=k\mp\sqrt{k^{2}-k} ,        \\
    \zeta^{\Delta} &=& 0 & \Leftrightarrow  &\mu &=&  0,    \\
    \eta^{\Delta} &=& 0 & \Leftrightarrow &\mu &=&0,k.
  \end{array}
\end{equation}
Thus the $\mu-$exceptional cases are $0, 1, k,
\overline{\mu}_{\pm}$ (compare with \cite{Ito01}).
  \item[$\underline{m=2}$] Note that $k \ge 1$
  \begin{equation}\label{eq:zeta eta HD exceptional m=2}
  \begin{array}{lcccccl}
    \zeta^{H} &=& 0 &\textrm{ never }, &&&         \\
    \eta^{H} &=& 0 & \Leftrightarrow &\mu &=& \displaystyle
\frac{1}{2},    \\
    \zeta^{\Delta} &=& 0 & \Leftrightarrow  &\mu &=&  0 ,   \\
    \eta^{\Delta} &=& 0 & \Leftrightarrow &\mu &=&0.
  \end{array}
\end{equation}
Thus the $\mu-$exceptional cases are $0, 1, \displaystyle
\frac{1}{2}$.
\end{description}
 Hence like for \textsc{Table \ref{tab:ricci mu exceptional cases m
ge 3}}, we can establish the \textsc{Table \ref{tab:ricci mu
exceptional cases m 1,2}}.
\end{rem}


\begin{small}
\begin{table}
\caption{\small Einstein equations, ~$m = 1, 2$,
~$\mu-$exceptional cases in
\textbf{Theorem \ref{thm:ricci sbcwp m ge 3}} }
\label{tab:ricci mu exceptional cases m 1,2}
\begin{tabular}{ccccccccc}
 \hline\hline
   $\mu$ & $m$ & $k$ & $\zeta^{H}$ & $\eta^{H}$ & $\zeta^{\Delta}$ &
$\eta^{\Delta}$
   & \begin{minipage}{.9 cm}
    \textit{ \centerline{genuine} \\\centerline{system}}
    \end{minipage} &  \begin{minipage}{.9 cm}
    \textit{ \centerline{formal} \\\centerline{system}}
    \end{minipage}
\\\\
 \hline
\\
  0 & $1$ & $\ge 1 $ & $-k$ & $-k$ & $0$ & $0$ &\eqref{eq:ricci
special 2} & \eqref{eq:ricci special 4}
\\\\
  1 & $1$ & $\ge 1 $ & $-[-1+k]$ & $-1+k$ & $1$ & $-1+k$ &
\eqref{eq:ricci special 2}& \eqref{eq:ricci special 4}
\\\\
  $k$& $1$ & $> 1 $ & 0 & $k(k-1)$ & $k$ & 0 & \eqref{eq:ricci
special 4} & -
\\\\
  $\overline{\mu}_{\pm}$& $1$ & $> 1 $ & $\displaystyle k
\frac{\overline{\mu}_{\pm}-1}{\overline{\mu}_{\pm}}$ & 0 &
$\overline{\mu}_{\pm}$ & $-k(\overline{\mu}_{\pm}-1)$
  & \eqref{eq:ricci special 4} & -
\\\\
  0 & 2 & $\ge 1 $  & $-k$ & $-k$  & 0 & 0 & \eqref{eq:ricci special
2} & \eqref{eq:ricci special 4}
\\\\
  1 & 2 & $\ge 1 $   & $-k$  & $k$  & 1 & $k$ & \eqref{eq:ricci
special 2}& \eqref{eq:ricci special 4}
\\\\
  $\displaystyle \frac{1}{2}$ & 2 & $\ge 1 $ & $-k$ & 0 & $\displaystyle\frac{1}{2}$ &
$\displaystyle\frac{k}{2}$ &\eqref{eq:ricci special 4} & -
\\\\
 \hline\hline
\end{tabular}
\end{table}
\end{small}


\subsection{Scalar Curvature} \label{subsec: sbcwp scalar curvature}

\begin{thm}
\label{thm:scalar sbcwp m ge 2} Let $B = (B_m,g_{B})$ and $F =
(F_k,g_{F})$ be two pseudo-Rieman\-nian manifolds with $m \geq 2$
and $k \geq 1$, $\displaystyle \mu \in \mathbb{R}\setminus
\left\{0, 1, \displaystyle-\frac{k}{m-1}
\right\}$ and $\psi \in C^\infty_{>0}(B)$. Then, the scalar
curvature $S$ of the base conformal warped product
\twistpar{B}{F}{\psi^\mu}{\psi}
verifies the relation 
\begin{equation}
   -\beta \Delta_{B}u + S_{B} u
   =
  S u^{2 \mu \alpha + 1}   - S_{F}u^{2(\mu - 1)\alpha + 1},
        \label{eq:scalar curvature m ge 2}
\end{equation}
where
\begin{equation}\label{eq:alpha m ge 3}
 \alpha =
\frac{2[k+(m-1)\mu]}{\{[k+(m-1)\mu]+(1-\mu)\}k + (m-2)\mu
[k+(m-1)\mu]},
\end{equation}
\begin{equation}\label{eq:beta m ge 3}
  \beta = \alpha 2[k+(m-1)\mu]
\end{equation}
and $\psi = u^{\alpha}>0$.

\end{thm}

\begin{proof}
 Applying \textbf{Theorem \ref{sca-c}} with $c=\psi^{\mu}$ and
$w=\psi$, we obtain
\begin{equation}\label{eq:scalar curv sbcwp}
\begin{array}{lcl}
\psi^{2\mu}S & = & S_B + S_F\psi^{2(\mu-1)} -\left[\displaystyle
2(m-1)\frac{\Delta_{B}\psi^{\mu}}{\psi^{\mu}}+ 2k
\frac{\Delta_{B} \psi}{\psi}\right]\\
& - & [(m-4)(m-1)   \mu^{2}
%
 +  2k(m-2) \mu
 +  k(k-1)]\displaystyle
\frac{g_B(\nabla^B\psi,\nabla^B\psi)}{\psi^{2}}.
\end{array}
\end{equation}
So by \eqref{eq:lappotencia}, with $t=\mu\neq 0, 1$, there results
\begin{eqnarray*}
\psi^{2\mu}S & = & S_B + S_F\psi^{2(\mu-1)}\\
& &-\left\{\left[2(m-1)+ \frac{\varsigma}{\mu(\mu-1)}
\right]\frac{\Delta_{B}\psi^{\mu}}{\psi^{\mu}}+
\left[2k-\frac{\varsigma}{\mu-1}\right] \frac{\Delta_{B}
\psi}{\psi}\right\},
\end{eqnarray*}
where $\varsigma = (m-4)(m-1)\mu^{2} +  2k(m-2) \mu +  k(k-1)$.
Hence, by \textbf{Lemma \ref{lem:2.1}} and \textit{Remark}
\textit{\ref{rem:operator lemma}} with
\begin{equation*}
  \begin{array}{lclcl}
    r_{1} &=& \displaystyle 2(m-1)+ \frac{\varsigma}{\mu(\mu-1)}, \\
    r_{2} &=& \displaystyle 2k-\frac{\varsigma}{\mu-1}
  \end{array}
\end{equation*}
and
\begin{equation}\label{eq:zeta eta particular}
  \begin{array}{lclcl}
    \zeta &=& r_{1} \mu+ r_{2} &=&2[k+(m-1)\mu] ,     \\
    \eta &=& r_{1} \mu^{2}+ r_{2}&=&\{[k+(m-1)\mu]+(1-\mu)\}k \\
    &&&&+
(m-2)\mu
[k+(m-1)\mu]\\
& & &=&\displaystyle \left\{\frac{\zeta}{2}+(1-\mu)\right\}k +
(m-2)\mu \frac{\zeta}{2},
  \end{array}
\end{equation}
we find
\begin{equation}\label{eq:scalar curv zeta-eta}
\begin{array} {lcl}
\psi^{2\mu}S & = & S_B + S_F\psi^{2(\mu-1)}\\
& & -\left[ \displaystyle (\eta - \zeta ) \frac{g_B(\nabla^B \psi
,\nabla^B \psi)}{\psi^{2}} + \zeta\frac{\Delta_{B}\psi}{\psi}
\right].
\end{array}
\end{equation}
Notice that (see also \eqref{eq:alpha condition} in the
\textbf{Appendix \ref{sec:appendix}})
\begin{equation}\label{eq:scalar eta}
 \eta = (m-1)(m-2)\mu^{2} + 2(m-2)k\mu +(k+1)k > 0 \textrm{ for
all } \mu \in \mathbb{R}.
\end{equation}
On the other hand, $\zeta = 0$ if and only if $\displaystyle \mu =
-\frac{k}{m-1}.$ Then the thesis follows by \textbf{Lemma
\ref{lem:2.1}} and taking
\begin{equation}\label{zeta eta scalar curv}
\displaystyle \alpha = \frac{\zeta}{\eta} \textrm{ and }  \beta =
\alpha \zeta.
\end{equation}
\end{proof}


The \textsc{Table 3}
is a synthesis of the cases not included in the \textbf{Theorem
\ref{thm:scalar sbcwp m ge 2}}. In that table $\zeta$ and $\eta$ are
computed with the expressions \eqref{eq:zeta eta particular} instead
of the originals in \textit{Remark} \textit{\ref{rem:operator
lemma}}. As above, this is the reason to include the column titled
``formal equation", and hence the equations written in that column
are justified a posteriori.

\medskip
\begin{small}
\begin{table}
 \addtocounter{table}{1}
\rotatebox[origin = c]{90}{
%
\begin{tabular}
{ccccccccc} 
  \multicolumn{9}{c}{
  \begin{minipage}{15 cm}
      \textsc{ \\ \small Table \arabic{table}.}
      {\textrm\small Scalar curvature equation,~$\mu-$exceptional
cases in
      \textbf{Theorem \ref{thm:scalar sbcwp m ge 2}}}
  \end{minipage}
}
\\
\\
  \hline\hline
  \\
$\mu$ & $m$ & $k$ & $\zeta$ & $\eta$ &
$
\begin{array}{ll}
    genuine \\
    equation\\
\end{array}%
$
%
%
& $
\begin{array}{ll}
    formal \\
    equation\\
\end{array}%
$
%
%
%
&$
\begin{array}{ll}
    equivalent \\
    equation\\
\end{array}%
$
%
%
%
&$
\begin{array}{ll}
    gometrical \\
    meaning\\
\end{array}%
$
%
%
  \\\\
  \hline
  \\
  0& $\ge 1$ & $0$ & $0$ & 0 & \eqref{eq:scalar curv sbcwp} &
\eqref{eq:scalar curv zeta-eta}
  &$S=S_{B}$&-
  \\
  0 & $\ge 1$ & $\ge 1$ & $2k$ & $(k+1)k$ & \eqref{eq:scalar curv
sbcwp}=\eqref{eq:scalar curv wp-1}
  &\eqref{eq:scalar curvature m ge  2}
  &\eqref{eq:scalar curv wp} & \textit{singly warped}
  \\
  1 & 1 & 0 & 0 & 0 &\eqref{eq:scalar curv sbcwp} & \eqref{eq:scalar
curv zeta-eta}
  & $S \equiv0$& -
  \\
  1 & 1 & 1 & 2 & 0 & \eqref{eq:scalar curv sbcwp} & \eqref{eq:scalar
curv zeta-eta}
  & \eqref{eq:conf 4},  $r=2$
  & \textit{conformal product}
  \\
  1 & 1 & $\ge 2$ & $2k$ & $k(k-1)$ & \eqref{eq:scalar curv sbcwp} &
\eqref{eq:scalar curv zeta-eta} &
   \eqref{eq:scalar curvature m ge  2}=\eqref{eq:yamabe usual
product}& \textit{conformal product}
  \\
  1 & 2 & 0 & 2 & 0 & \eqref{eq:scalar curv sbcwp} & \eqref{eq:scalar
curv zeta-eta} & \eqref{eq:conf 4}, $r=2$
  &
$\begin{cases}
    \textit{base conformal}\\
    \textit{Nirenberg pb. type}
\end{cases}$%
  \\
  1 & $\ge 3$ & $0$ & $2(m-1)$ & $(m-1)(m-2)$ & \eqref{eq:scalar curv
sbcwp} & \eqref{eq:scalar curv zeta-eta}
  & \eqref{eq:conf 1}, $r=2$&
$\begin{cases}
  \textit{base conformal}\\
  \textit{Yamabe eq. type}
\end{cases}$
  \\
  1 & $\ge 2$ & $\ge 1$ & $2[k+m-1]$ & $\displaystyle
\frac{\zeta}{2}\left(\frac{\zeta}{2}-1\right)$ & \eqref{eq:scalar
curv sbcwp}
  & \eqref{eq:scalar curv zeta-eta} & \eqref{eq:scalar curvature m ge
 2}=\eqref{eq:yamabe usual product}
  &\textit{conformal product}
  \\
  $\displaystyle -\frac{k}{m-1}$ & $\ge 2$ & $\ge 1$ & 0 & $>0$ &
  \eqref{eq:scalar curv zeta-eta} & - & \eqref{eq:rem scalar
curvature  -k/(m-1)}& -
  \\\\
  $\neq \displaystyle \frac{k+1}{2}, 0, 1$ & $1$ & $\ge 1$ & $2k$ &
$k(k + 1 - 2 \mu)$
  & \eqref{eq:scalar curvature m ge  2}
   & - &  - & -
  \\ \\
  $\displaystyle \frac{k+1}{2}$ & $1$ & $> 1$ & $2k$ & $0$
  & \eqref{eq:scalar curv zeta-eta}
   & - &  - & -
  \\\\
\hline\hline
\end{tabular}
%
}
%
%
\label{tab:scalar curvature mu exceptional cases m ge 1}
\end{table}
\end{small}
%
%
%


%
%
\bigskip

All the other cases are covered in \textbf{Theorem \ref{thm:scalar
sbcwp m ge 2}}.


\begin{rem}\label{rem:mu special} We want to make some comments about
the results in the \textsc{Table 3} where we have three important
cases:
\begin{description}
  \item[($\mu = 0$)] As it was mentioned in \S 1,
this case corresponds exactly to standard warped products. The
relation \eqref{eq:scalar curvature m ge 2} is well defined and
reproduced in \eqref{eq:scalar curv wp}.

  \item[($\mu = 1, k = 0, m \geq 3 $)] This situation corresponds to a
  conformal change in the base. Again \eqref{eq:scalar curvature m ge
2} is well defined and now reproduces \eqref{eq:conf 1} with
$r=2$, and hence \eqref{eq:yamabe} too.

  \item[($\mu = 1, k,m \geq 1, k + m \geq 3$)](i.e., \textit{rows 5 or 8})
  We have a conformal change in the usual
  product, more explicitly, $(B \times F,g = \psi^{2} (g_{B} +
g_{F})) $. In this case \eqref{eq:scalar curvature m ge 2} is well
defined also, and reproduce with $\displaystyle \alpha =
\frac{2}{m + k -2}$ and $ \displaystyle \beta = 4 \frac{m + k -1}{
m + k -2}$, the equation
\begin{equation}
    - 4 \frac{m + k -1}{
m + k -2} \Delta_{g_{B}} u + (S_{g_{B}} +S_{g_{F}} )u
    = S u^{1+\frac{4}{m+k-2}},
    \label{eq:yamabe usual product}
\end{equation}
where $g = u^{\frac{4}{m+k-2}} (g_{B} + g_{F})$, $u \in
C^\infty_{>0}(B)$, $\psi = u^{\frac{2}{m+k-2}}$ and
$c_{m+k}=\beta$.
\end{description}

Now we will analyze the cases included neither in the previous
items nor in \textbf{Theorem \ref{thm:scalar sbcwp m ge 2}}.

\begin{description}
  \item[$(m = 1)$] $ $Let $k \geq 1$. It is clear that the
involved differential equations are ordinary and $S_{B}\equiv 0$. If
\begin{itemize}
  \item $\displaystyle\left(\mu \neq 0, 1, \frac{k+1}{2}\right)$
  By the same proof of \textbf{Theorem \ref{thm:scalar sbcwp m ge 2}},
  the equation \eqref{eq:scalar curvature m ge  2}
  is valid.
  \item $\displaystyle(\mu = 1, k \geq 2)$ It is a particular case
  of the above item ($\mu = 1, k,m \geq 1, k + m \geq 3$), so
\eqref{eq:scalar curvature m ge  2}
   is true again.
  \item $\displaystyle\left(\mu = \frac{k+1}{2}, k \neq 1\right)$ It
  is
  possible to apply \eqref{eq:scalar curv zeta-eta} so
\begin{equation}
 \label{eq: rem (k+1)/2, k >1}
  \psi^{k+1}S = 2k\left(-\frac{\Delta_{B} \psi }{\psi} +
\frac{|\nabla^{B}
  \psi|_{B} ^{2}}{\psi^{2}} \right) + S_{F} \psi^{k-1}
\end{equation}
  \item $\displaystyle\left(\mu = \frac{k+1}{2}, k = 1\right)$
Clearly $\mu =
  1$, hence \eqref{eq:scalar curv sbcwp} results by applying
\eqref{eq:conf 4} with $r =
  2$, i.e.
\begin{equation}\label{eq: rem m = 1, k = 1, mu = 1}
  \psi^{2} S = 2 \left(-\frac{\Delta_{B} \psi}{\psi} +
\frac{|\nabla^{B}
  \psi|_{B} ^{2}}{\psi^{2}}\right).
\end{equation}
Confront with the precedent case.

\end{itemize}

  \item[$\left(\displaystyle m \geq 2, \mu = -\frac{k}{m-1}\right)$]
In this
  case by \eqref{eq:scalar curv zeta-eta}  the relation among the
scalar curvatures is
%
%
\begin{equation}\label{eq:rem scalar curvature -k/(m-1)}
  -k\left[1 + \frac{k}{m-1}\right] \frac{|\nabla^{B}
  \psi|_{B} ^{2}}{\psi^{2}}= \psi^{-2\frac{k}{m-1}} S - S_{B} -
    S_{F} \psi^{-2\left(1 + \frac{k}{m-1}\right)}.
\end{equation}
\end{description}
\end{rem}

\begin{rem} \label{rem:beta > 0} Note that $\beta > 0$ in
\textbf{Theorem \ref{thm:scalar sbcwp m ge 2}}, while this is not
always true if $m=1$.
\end{rem}

\begin{proof}
\textbf{\textbf{(of Theorem \ref{thm:scurv conf warped m ge 2})}}
It is an immediate consequence of the above results of this
section.
\end{proof}

\section{The nonlinearities in a \bcwpar{\psi,\mu}
scalar curvature relations.
} \label{sec:constant scalar curvature bcwp with m ge 2}

In this section, we will mainly consider some general properties
of the nonlinear partial differential equation in \eqref{eq:scalar
curvature m ge 2}, regarding especially the type of
nonlinearities. The main aim of this study is to deal with the
question of existence and multiplicity of solutions for problem
\textbf{(Pb-sc)}.
The corresponding results will be presented in forthcoming
articles (see \cite{DobarroUnal07}).

From now on, we will denote by $\discr{\cdot} $, the discriminant
of a quadratic polynomial in one variable.

\subsection{Base $B_m$ with dimension $m\ge 2$}
\label{subsec: General Properties}

\begin{rem}
\label{rem: classification of nonlinerities m ge 2} Under the
hypothesis of \textbf{Theorem \ref{thm:scalar sbcwp m ge 2}}. In
order to classify the type of non linearities involved in
\eqref{eq:scalar curvature m ge 2},
we will analyze the exponents as a function of the parameter $\mu$
and the dimensions of the base $m \ge 2$ and  of the fiber $k \ge
1$
(see \textsc{Table \ref{tab:scalar curvature nonlinearity type m
ge 2}} below).

\medskip

Note that by \eqref{eq:scalar eta},
$\alpha > 0$ if and only if $\displaystyle \mu > -\frac{k}{m-1}$
and by the hypothesis $\displaystyle \mu \neq -\frac{k}{m-1}$ in
\textbf{Theorem \ref{thm:scalar sbcwp m ge 2}}
, results $\alpha \neq 0$.

\medskip

We now introduce the following notation:
\begin{equation}\label{eq:p,q}
  \begin{array}{ll}
    p= p(m,k,\mu)= & 2 \mu \alpha + 1 \textrm{ and }\\
    q= q(m,k,\mu)= & 2(\mu - 1)\alpha + 1 = p - 2 \alpha,
  \end{array}
\end{equation}
where $\alpha$ is defined by \eqref{eq:alpha m ge 3}.

\medskip

Thus, for all $ m, k, \mu $ as above, $p > 0$. Indeed, by
\eqref{eq:scalar eta}, $p > 0$ if and only if $\varpi
> 0$, where
\begin{equation*}\label{}
 \begin{array}{ll}
   \varpi &:= \varpi (m,k,\mu) \\
      &:=4 \mu [k + (m-1)\mu] + (m-1)(m-2)\mu^{2} + 2(m-2)k\mu
+(k+1)k \\
          &=(m-1)(m+2)\mu^{2} + 2 m k \mu  + (k+1)k.
 \end{array}
\end{equation*}
But $\discr{\varpi} \le -4km^{2} \le -16$ and $m>1$, so $\varpi
>0.$

\medskip

Unlike $p$, $q$ changes sign depending on $m$ and $k$.
Furthermore, it is important to determine the position of $p$ and
$q$ with respect to $1$ as a function of $m$ and $k$. In order to
do that, we define
\begin{equation}\label{eq:m,k condition}
  D:=\{(m,k)\in \mathbb{N}_{\ge 2}\times\mathbb{N}_{\ge 1}:
\discr{\varrho(m,k,\cdot)}<0\},
\end{equation}
where $\mathbb{N}_{\ge l} := \{j \in \mathbb{N}: j \ge l\}$,
\begin{equation*}\label{}
 \begin{array}{ll}
   \varrho &:= \varrho (m,k,\mu) \\
       &:=4 (\mu-1) [k + (m-1)\mu] + (m-1)(m-2)\mu^{2} +
2(m-2)k\mu +(k+1)k \\
           &=(m-1)(m+2)\mu^{2} + 2 (m k -2(m-1))\mu  + (k-3)k
 \end{array}
\end{equation*}
and the discriminant of $\varrho (m,k,\cdot)$ is
\begin{equation*}\label{eq:discriminant}
    \discr{\varrho(m,k,\cdot)}=-4((m-2)k - 4(m-1))(k+m-1).
\end{equation*}
Note that by \eqref{eq:scalar eta}, $q>0$ if and only if $\varrho
>0$. Furthermore $q=0$ if and only if $\varrho
=0$. But here $\discr{\varrho(m,k,\cdot)}$ changes its sign as a
function of $m$ and $k$.


In \textsc{Table \ref{tab:scalar curvature nonlinearity type m ge
2}} below, we denote $\mathcal{C}D = (\mathbb{N}_{\ge 2} \times
\mathbb{N}_{\ge 1}) \setminus D$ if $D \subseteq \mathbb{N}_{\ge
2} \times \mathbb{N}_{\ge 1}$ and $\mathcal{C}I = \mathbb{R}
\setminus I$ if $I \subseteq \mathbb{R}$. If $(m,k) \in
\mathcal{C}D$, let $\mu_{-}$ and $\mu_{+}$ the two (eventually
one, see \textit{Remark \ref{rem:m=6 and k=5}} below) roots of
$q$, $\mu_{-} \le \mu_{+}$. Besides, if
$\discr{\varrho(m,k,\cdot)}>0 $, then $\mu_{-}<0$; unlike
$\mu_{+}$ can take any sign.


\begin{small}
\begin{table}
\caption{\small
Nonlinearities in scalar curvature equation type \eqref{eq:scalar
curvature m ge 2} for $m \ge 2$, see \textit{Notation \ref{nota:
right columns in tables 45678} } }
 \label{tab:scalar curvature nonlinearity type m ge 2}
\begin{tabular}{ccccc}
\hline \hline
   $(m,k)\in $
   &  $\mu \in $
   &$ \alpha $
   & $p,q$
   &
    %
    %
    \begin{minipage}{1.5 cm}
    \textit{ \centerline{type of } \\ \centerline{$p,q$ non- }  \\
    \centerline{linearity}
    }
    \end{minipage}
   \\
 \hline \\
   $\mathbb{N}_{\ge 2}\times\mathbb{N}_{\ge 1}$
   & $\left(-\infty,-\displaystyle\frac{k}{m-1}\right)$
   & $<0$
   &$1<p<q$
   & \textit{super-lin}
   \\\\
   $D$
   & $\left(-\displaystyle\frac{k}{m-1},0\right)$
   & $0<$
   &$0<q<p<1$
   & \textit{sub-lin}
   \\\\
  $\mathcal{C}D $
  &  $\left(-\displaystyle\frac{k}{m-1},0\right)\cap (\mu_{-} ,
\mu_{+})$
  & $0<$
  &$q<0<p<1$
  & $\left\{\begin{minipage}{1.5 cm}
   \textit{ \centerline{sub-lin} \\ \centerline{sing}}
   \end{minipage}\right.$
  \\\\
  $\mathcal{C}D $
  & $\left(-\displaystyle\frac{k}{m-1},0\right)\cap
\mathcal{C}[\mu_{-} , \mu_{+}]$
  & $0<$
  &$0<q<p<1$
  &
   \textit{sub-lin}
  \\\\
  $\mathcal{C}D $
  & $\left(-\displaystyle\frac{k}{m-1},0\right)\cap \{\mu_{-} ,
\mu_{+}\}$
  & $0<$
  & $q=0<p=2\alpha<1$
  &
   $\left\{\begin{minipage}{1.5 cm}
   \textit{ sub-lin \\ \centerline{\underline{non}-hom}}
   \end{minipage}\right.$
  \\\\
  $D $
  & $(0,1)$
  & $0<$
  & $0<q<1<p$
  &
   $\left\{\begin{minipage}{1.5 cm}
   \textit{ \centerline{super-lin} \\ \centerline{sub-lin}}
   \end{minipage}\right.$
  \\\\
  $\mathcal{C}D $
  & $(0,1)\cap (\mu_{-} , \mu_{+})$
  & $0<$
  & $q<0<1<p$
  &
   $\left\{\begin{minipage}{1.5 cm}
   \textit{ \centerline{super-lin} \\ \centerline{sing}}
   \end{minipage}\right.$
  \\\\
  $\mathcal{C}D $
  & $(0,1)\cap \mathcal{C}[\mu_{-} , \mu_{+}]$
  & $0<$
  & $0<q<1<p$
  &
   $\left\{\begin{minipage}{1.5 cm}
   \textit{ \centerline{super-lin} \\ \centerline{sub-lin}}
   \end{minipage}\right.$
  \\\\
  $\mathcal{C}D $
  &  $(0,1)\cap \{\mu_{-} , \mu_{+}\}$
  & $0<$
  & $q=0<1<p=2\alpha$
  &
   $\left\{\begin{minipage}{1.5 cm}
   \textit{\centerline{super-lin} \\ \centerline{\underline{non}-hom}}
   \end{minipage}\right.$
  \\\\
  $\mathbb{N}_{\ge 2}\times\mathbb{N}_{\ge 1}$
  & $(1,+\infty)$
  & $0<$
  & $1<q<p$
  & \textit{super-lin}
  \\\\
  \hline\hline
\end{tabular}
\end{table}
\end{small}



We remark that all the rows in \textsc{Table \ref{tab:scalar
curvature nonlinearity type m ge 2}} are \textit{nonempty}, this
means that the conditions established in each row are verified for a
suitable choice of the parameters and manifolds. On the other hand,
we observe that $\beta$ is always positive as it was mentioned in
\textit{Remark \ref{rem:beta > 0}}.  Note that for any row in
\textsc{Table \ref{tab:scalar curvature nonlinearity type m ge 2}},
the corresponding type of nonlinearity suggested by the exponents is
modified by the scalar curvature of the fiber, $S_{F}$ and by the
function $S$.

Furthermore, depending on whether the base is Riemannian or not,
then the linear part is elliptic or not, respectively.

\end{rem}

\begin{nota}
\label{nota: right columns in tables 45678} In the last right hand
side columns of \textsc{Tables \ref{tab:scalar curvature
nonlinearity type m ge 2}}, \ref{tab:scalar curvature nonlinearity
type m = 1, k ge 4}, \ref{tab:scalar curvature nonlinearity type m
= 1, k = 3}, \ref{tab:scalar curvature nonlinearity type m = 1, k
= 2} and \ref{tab:scalar curvature nonlinearity type m = 1, k =
1}, we will use the notation explained below:
\begin{itemize}
  \item \textit{super-lin} means that the corresponding exponent $ > 1$,
roughly speaking
  \textit{super-linear}
  \item \textit{sub-lin} means that the corresponding exponent $> 0$
  and
  $< 1$, roughly speaking
  \textit{sub-linear}
  \item \textit{\underline{non}-hom} means that the corresponding exponent
  $=0$, roughly speaking \textit{non-homogeneous}
  \item \textit{sing} means that the corresponding exponent $ < 0$, roughly
  speaking \textit{singular}.
\end{itemize}
However, all these conditions depend strongly on the corresponding
coefficients in the whole specific non-linearity. More clearly,
we can say that
the right columns of the tables mentioned above are exact when $S$
and $S_F$ are strictly positive constants.
\end{nota}

\begin{rem}\label{rem:m=6 and k=5}
Note that when we consider $\discr{\varrho}=0$, we look for
solutions $(m,k)\in \mathbb{N}_{\ge 2} \times \mathbb{N}_{\ge 1}$,
in particular ordered pairs with natural components. It is easy to
see that
\begin{equation}\label{eq:m,k condition-level 0}
  \begin{array}{ccl}
    D_{0} & = &  \{(m,k)\in \mathbb{N}_{\ge 2}\times\mathbb{N}_{\ge
1}:
  \discr{\varrho(m,k)}=0\}\\
     & = & \left\{(m,k)\in \mathbb{N}_{\ge 3}\times\mathbb{N}_{\ge
1}:  k=4\displaystyle
     \frac{m-1}{m-2}\right\}\\
     & = & \{(3,8),(4,6),(6,5)\}.
  \end{array}
\end{equation}
All the other solutions of $\discr{\varrho}=0$ in $\mathbb{R}^2$
have no natural components.

Then, for $(m,k)=(3,8) \in D_0$ $((4,6), (6,5)\textrm{
respectively })$, $\displaystyle -\frac{k}{m-1}$ takes the value
$-4$ $(-2, -1 \textrm{ respectively })$ and $\mu_{-}=\mu_{+}=-2$
$\left( \displaystyle -1,-\frac{1}{2}\right.$ respectively).
%
%
In such a case, when $\mu=\mu_{-}=\mu_{+}$, the fifth row in
\textsc{Table \ref{tab:scalar curvature nonlinearity type m ge 2}}
establishes that $q=0$, $p=\displaystyle \frac{1}{3}$ $\left(
\displaystyle \frac{1}{2},\frac{2}{3}\textrm{ respectively}\right)$,
$\alpha =\displaystyle \frac{1}{6}$ $\left(\displaystyle
\frac{1}{4},\frac{1}{3} \textrm{ respectively}\right)$ and $\beta
=\displaystyle \frac{4}{3}$ $\left( \displaystyle
\frac{3}{2},\frac{10}{3}\textrm{ respectively}\right)$.
%

Note that for the elements in $D_{0}$, the sum of the two
components is either $11$ or $10$, both particularly interesting
values in the physical applications. More precisely in the
problems of the extra dimensions in cosmology, super-gravity and
string theory (i.e. see \cite{Aharony-Gubser-Maldacena-Ooguri-Oz,
Argurio98, Gauntlett-Kim-Waldram01, Gauntlett-Kim-Waldram01-b,
Gauntlett-Kim-Pakis-Waldram02,
Gauntlett-Martelli-Sparks-Waldram04, Randall-Sundrum99a,
Randall-Sundrum99b}).
\end{rem}

\begin{nota}\label{nota: yamabe exponent}
From now on, for $m \ge 3$ we will denote the Sobolev critical
exponent by $2^*=\displaystyle{\frac{2m}{m-2}}$ and
$\displaystyle{p_Y=q_Y= \frac{4}{m-2}+1=\frac{m+2}{m-2}}=2^* - 1$.
\end{nota}

\begin{rem}\label{rem:mu and sobolev}
Let $m \ge 3$. Now we will show that there exist particular values
$\mu_{p_Y}$ and $\mu_{q_Y}$ such that the position of $\mu$ with
respect to them, indicates that the corresponding $p$ or $q$ are
sub-critical, critical or super-critical. The critical and
super-critical cases will correspond to the conditions in the
first row of \textsc{Table \ref{tab:scalar curvature nonlinearity
type m ge 2}}. Indeed, by an easy but lengthy computation we have
\begin{description}
  \item[$\underline{p>p_Y}$]if and only if $\mu < \mu_{p_Y}=
  \displaystyle{-\frac{k+1}{m-2}}$.
  \item[$\underline{q>q_Y}$]if and only if $\mu <\mu_{q_Y}=
  \displaystyle{-\frac{k}{m-2}}$.
\end{description}
Moreover,
\begin{description}
  \item[$\underline{p=p_Y}$] is verified if and only if $\mu =
  \displaystyle{-\frac{k+1}{m-2}}$; and consequently $\alpha =
-\displaystyle{\frac{2}{k+1}}$,
   $\beta = \displaystyle{4\frac{m-1}{m-2}-4\frac{k}{k+1}}>0$ and
  $q = \displaystyle{p_Y+\frac{4}{k+1}}$. Hence the
  equation \eqref{eq:scalar curvature m ge 2} takes the form
  \begin{equation}
   \displaystyle{-\Big(4\frac{m-1}{m-2}-4\frac{k}{k+1}\Big)}
\Delta_{B}u + S_{B} u
   =
  S u^{p_Y}   - S_{F}u^{p_Y+ \frac{4}{k+1}}.
        \label{eq:scalar curvature m ge 3 pY}
\end{equation}
  \item[$\underline{q=q_Y}$] is verified if and only if $\mu = \mu_{q_Y}
  =
  \displaystyle{-\frac{k}{m-2}}$; and consequently \\ $\alpha =
  -\displaystyle{\frac{2}{k+m-2}}$,
  $\beta = \displaystyle{\frac{4k}{(k+m-2)(m-2)}}>0$ and \\
  $p = q_Y \displaystyle{- \frac{2}{k+m-2} }$. Hence the
  equation \eqref{eq:scalar curvature m ge 2} takes the form
  \begin{equation}
   \displaystyle{-\frac{4k}{(k+m-2)(m-2)}} \Delta_{B}u + S_{B} u
   =
  S u^{q_Y - \frac{2}{k+m-2}}   - S_{F}u^{q_Y}.
        \label{eq:scalar curvature m ge 3 qY}
\end{equation}
Note that $\mu_{q_Y}$ is the exceptional value $\overline{\mu}$ in
\textbf{Theorem \ref{thm:ricci sbcwp m ge 3}} (see \textsc{Table
\ref{tab:ricci mu exceptional cases m ge 3}}).
\end{description}

We observe also that $\mu_{p_Y}<\mu_{q_Y}<\displaystyle
-\frac{k}{m-1}$, so that at least one of the two exponents is no
sub-critical only if we stay in the conditions of the first row in
the \textsc{Table \ref{tab:scalar curvature nonlinearity type m ge
2}}.
\end{rem}

\begin{rem}\label{rem:mu pm-infty and sobolev}
Let $m \ge 3$. Now, we will study the behavior of Equation
\eqref{eq:scalar curvature m ge 2}, when $\mu \longrightarrow
\pm\infty$. Consider $\mu \longrightarrow \pm \infty$ , then by
\eqref{eq:alpha m ge 3} we have (see \textsc{table \ref{tab:scalar
curvature nonlinearity type m ge 2}})
\begin{equation}\label{eq:alpha bis m ge 3}
 \alpha  =\displaystyle
\frac{2}{
\left\{\displaystyle
1+\frac{\displaystyle\frac{1}{\mu}-1}{\displaystyle\frac{k}{\mu}+
m - 1}\right\}
k + (m-2)\mu  }\longrightarrow \pm 0
\end{equation}
and
\begin{equation}\label{{eq:alpha.mu m ge 3}}
  \alpha\mu \longrightarrow \displaystyle \frac{2}{m-2}.
\end{equation}
Hence,
\begin{equation}\label{eq: going to yamabe}
  \begin{array}{lll}
    \beta & = &\alpha 2 [k + (m-1)\mu] =  \alpha 2 k + 2(m-1) \alpha
\mu  \longrightarrow \beta_{Y}=
    4\displaystyle\frac{m-1}{m-2}\\
    p     & = & 2 \mu \alpha + 1 \longrightarrow \displaystyle
p_{Y}=\frac{4}{m-2} + 1\\
    q     & = & 2 (\mu - 1) \alpha + 1 = p - 2 \alpha
    \longrightarrow \displaystyle q_{Y}=\frac{4}{m-2} + 1,
  \end{array}
\end{equation}
with $q_{Y}=p_{Y}$. Thus, roughly speaking the limit equation of
\eqref{eq:scalar curvature m ge 2} for $\mu \longrightarrow \pm
\infty$ results
\begin{equation}
   - 4\displaystyle\frac{m-1}{m-2}\Delta_{B}u + S_{B} u
   =
  (S - S_{F})u^{\frac{4}{m-2} + 1},
        \label{eq:going to yamabe-2}
\end{equation}
by ``a suitable definition of $S$". Notice the similarity of this
equation with the Yamabe type equation associated to a conformal
change in the base (see equation \eqref{eq:yamabe}). Furthermore,
by the last part of \textit{Remark \ref{rem:mu and sobolev}}, the
approximation is by super-critical problems when $\mu
\longrightarrow -\infty$ and by sub-critical problems when $\mu
\longrightarrow +\infty$.

\begin{figure}
   \epsfig{file=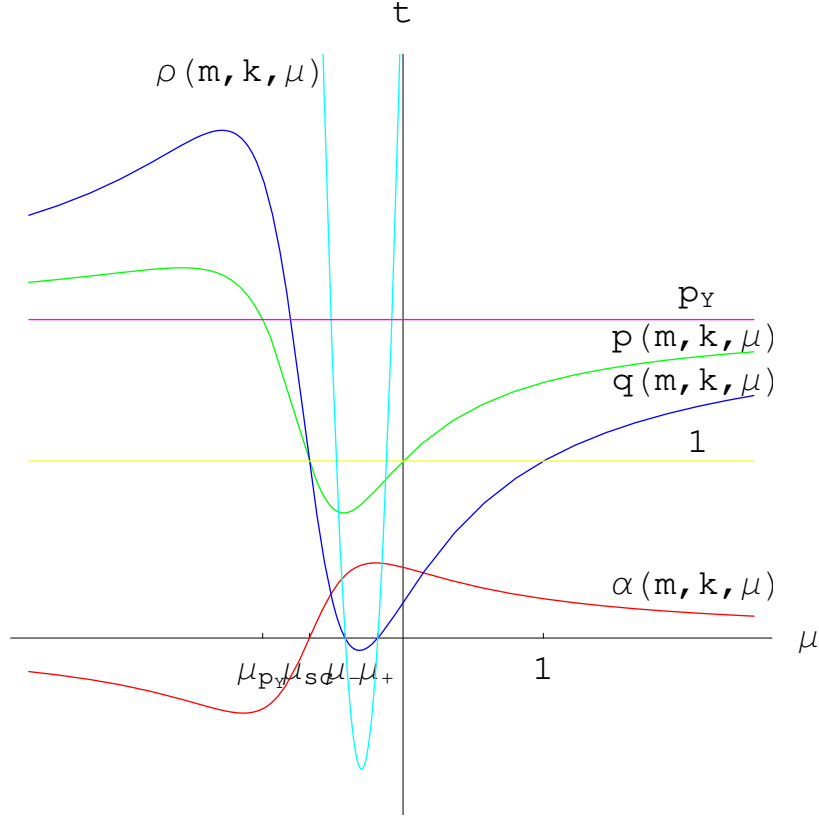,width=11cm,angle=0}
   \begin{small}
   \caption{Example:
   $(m,k)=(6,4) \in \mathcal{C}D$
   }
   \end{small}
\end{figure}

\end{rem}

\subsection{Base $B_m$ with dimension $m = 1$ }
\label{subsec:constant scalar curvature bcwp with m = 1}

\begin{rem}
\label{rem: classification of nonlinerities m = 1} As in the case
of \textit{Remark \ref{rem: classification of nonlinerities m ge
2}}, we will classify the type of non linearities involved in
\eqref{eq:scalar curvature m ge 2}, obviously when this equation
is verified \Big(see \textit{Remark \ref{rem:mu special}} and
cases either $k \ge 1$ and $\mu \neq 0, 1, \displaystyle
\frac{k+1}{2} $ or $k \ge 2$ and $\mu = 1$ there \Big).
Furthermore $m=1$ implies that the equations in \eqref{eq:scalar
curvature m ge 2} are ordinary differential equations and that the
curvature tensor of the base is $0$, and consequently $S_{B}\equiv
0$. Analogously, $S_{F}\equiv 0$ if $k=1$. Hence, we will analyze
the exponents as a function of the parameter $\mu$ and the
dimension of the fiber $k \ge 1$.

\noindent Similar to the case of $m \ge 2$, for any row in the
\textsc{Tables
\ref{tab:scalar curvature nonlinearity type m = 1, k ge 4},
\ref{tab:scalar curvature nonlinearity type m = 1, k = 3},
\ref{tab:scalar curvature nonlinearity type m = 1, k = 2},
\ref{tab:scalar curvature nonlinearity type m = 1, k = 1},
}the corresponding type of nonlinearity is modified by the scalar
curvature of the fiber $S_{F}$ and by the function $S$.

The problem \textbf{(Pb-sc)} for $m=1$ and the corresponding
nonlinear ordinary differential equations for low values of $k$ are
particularly interesting in physical applications (see
\cite{Ito01,Ito01-b,Ito03}, Kaluza-Klein theory and Randall-Sundrum
theory).

\medskip

By these hypothesis, we have
\begin{equation}\label{eq: alpha m=1}
 0 \neq \alpha = \frac{2}{-2\mu + k +1} = \frac{1}{-\mu + k_{1}}
\end{equation}
and
\begin{equation}\label{eq: beta m=1}
 0 \neq \beta = \frac{4k}{-2\mu + k +1}=\frac{2k}{-\mu + k_{1}},
\end{equation}
 where
\begin{equation}\label{eq:k1}
  1 \le k_{1}:= \frac{k+1}{2}.
\end{equation}
Note that by \eqref{eq: alpha m=1}, we have that $\alpha > 0$ if
and only if $\mu < k_{1}$. By \eqref{eq: beta m=1}, we also have
that $\beta >0$ if and only if $\mu < k_{1}$.

Furthermore by the same notation introduced in \eqref{eq:p,q}, we
have
\begin{equation}\label{eq:p m=1}
    p= p(1,k,\mu)=  2 \mu \alpha + 1
    =  \displaystyle\frac{\mu + \displaystyle\frac{k+1}{2}}{-\mu +
\displaystyle\frac{k+1}{2}}
    =  \displaystyle \frac{\mu + k_{1}}{-\mu + k_{1}}
\end{equation}
and
\begin{equation}\label{eq: q m=1}
  q= q(1,k,\mu)=  2(\mu - 1)\alpha + 1 = p - 2 \alpha
  = \displaystyle\frac{\mu + \displaystyle \frac{k-3}{2}}{-\mu +
\displaystyle\frac{k+1}{2}}
  =  \displaystyle \frac{\mu + k_{1}-2}{-\mu + k_{1}}.
\end{equation}
In particular,
\begin{description}
  \item[i] $\mu > k_{1}$ if and only if $\alpha < 0$ if and
only if $p < q$.
  \item[ii] $p < 1$ if and only if $\mu \alpha<
0$ and $q < 1$ if and only if $(\mu-1) \alpha< 0$.
  \item[iii] $p > 0$ if and only if $\mu \in (-k_{1},k_{1})$.
  \item[iv] $q>0$ if and only if $\mu \in (2-k_{1},k_{1})$ or $\mu
  \in (k_{1},2-k_{1})$.
  \item[v] $ 2 - k_{1}<0$ if and only if $3<k$.
\end{description}

\medskip
\end{rem}

Now we will separately analyze the cases $k \ge 4$, $k=3$, $k=2$
and  $k=1$ (see \textbf{v.} above and the first paragraph of this
subsection).


\begin{description}
  \item[$\underline{k \ge 4}$] then $\displaystyle 2-k_{1} \le
-\frac{1}{2} < 0 < \frac{5}{2} \le
  k_{1}$. Thus we obtain \textsc{Table
\ref{tab:scalar curvature nonlinearity type m = 1, k ge 4}.
}
\end{description}

\begin{small}
\begin{table}
\caption{\small Nonlinearities in scalar curvature equation type
\eqref{eq:scalar curvature m ge 2} for $m = 1 $ and $k \ge 4$, see
\textit{Notation \ref{nota: right columns in tables 45678}}  }
  \label{tab:scalar curvature nonlinearity type m = 1, k ge 4}
\begin{tabular}{cccc}
\hline \hline
   $\mu \in $ &$ \alpha \in $
   & $p,q$ &
    %
    %
   \begin{minipage}{1.5 cm}
   \textit{ \centerline{type of } \\\centerline{p,q non-}\\ \centerline{linearity}
   }
   \end{minipage}
   \\
 \hline
   \\
   $(-\infty,-k_{1})$ & $\Big(0,\displaystyle \frac{1}{2k_{1}}\Big)$
   &$q<p<0$  & \textit{sing}
   \\\\
   $ \{-k_{1}\}$ & $\Big\{\displaystyle \frac{1}{2k_{1}}\Big\}$
   &$q<p=0<1$  &
   $\left\{\begin{minipage}{1.5 cm}
   \textit{\centerline{\underline{non}-hom} \\ \centerline{sing}}
   \end{minipage}\right.$
   \\\\
   $ (-k_{1},2-k_{1})$
   & $\Big(\displaystyle \frac{1}{2k_{1}},\frac{1}{2(k_{1}-1)}\Big)$
& $q<0<p<1$  & 
  $\left\{\begin{minipage}{1.5 cm}
   \textit{\centerline{sub-lin} \\ \centerline{sing}}
   \end{minipage}\right.$
  \\\\
   $ \{2-k_{1}\}$ & $\Big\{\displaystyle \frac{1}{2(k_{1}-1)}\Big\}$
& $\displaystyle q=0<p=\frac{1}{k_{1}-1}<1$ &
%
   $\left\{   \begin{minipage}{1.5 cm}
   \textit{\centerline{sub-lin} \\ \centerline{\underline{non}-hom}}
   \end{minipage}\right.$
   \\\\
   $ (2-k_{1},0)$ & $\Big(\displaystyle
\frac{1}{2(k_{1}-1)},\frac{1}{k_{1}}\Big)$ & $0<q<p<1$ &
\textit{sub-lin}
  \\\\
   $(0,1)$ & $\Big(\displaystyle
\frac{1}{k_{1}},\frac{1}{k_{1}-1}\Big) $
   &
   $ 0<q<1<p$  &
   $\left\{   \begin{minipage}{1.5 cm}
   \textit{\centerline{super-lin} \\ \centerline{sub-lin}}
   \end{minipage}\right.$
   \\\\
  $\{1\}$ & $\displaystyle \frac{1}{k_{1}-1}$
  & $q=1<p=\displaystyle \frac{k_{1}+1}{k_{1}-1}$   &
  $\left\{   \begin{minipage}{1.5 cm}
  \textit{\centerline{super-lin} \\ \centerline{lin}}
  \end{minipage}\right.$
   \\\\
  $(1,k_{1})$ & $\Big(\displaystyle \frac{1}{k_{1}-1},+\infty\Big)$
  &$1<q<p $  & \textit{super-lin }
   \\\\
  $(k_{1},+\infty)$ & $(-\infty,0) $
  &$p<q<0 $   & \textit{ sing }
   \\\\
  \hline\hline
\end{tabular}
\end{table}
\end{small}

\begin{description}
  \item[$\underline{k = 3}$] this implies $\displaystyle 2-k_{1} = 0
< k_{1}=2$. Hence we have \textsc{Table \ref{tab:scalar curvature
nonlinearity type m = 1, k = 3}}.
\end{description}


\begin{small}
\begin{table}
  \caption{\small Nonlinearities in scalar curvature equation type
\eqref{eq:scalar curvature m ge 2} for $m = 1 $ and $k = 3$, see
\textit{Notation \ref{nota: right columns in tables 45678}} }
  \label{tab:scalar curvature nonlinearity type m = 1, k = 3}
\begin{tabular}{cccc}
\hline\hline
\\
   $\mu \in $ &$ \alpha \in$
   %
   %
   & $p,q$ &
    \textit{ type of $p,q$ non-linearity}
    %
   \\\\
 \hline
   \\
 $(-\infty,-2)$ & $\displaystyle \Big(0,\frac{1}{4}\Big)$
   %
   %
   &$q<p<0$  & \textit{sing}
   \\\\
 $ \{-2\}$ & $\displaystyle\Big\{\frac{1}{4}\Big\}$
   %
   %
   &$\displaystyle q=-\frac{1}{2}<p=0$  &
\textit{\underline{non}-hom / sing}
   \\\\
   $(-2,0)$ & $\displaystyle\Big(\frac{1}{4},\frac{1}{2}\Big)$
   %
   %
   &$q<0<p<1$  &
   \textit{sub-lin / sing}
   \\\\
   $(0,1)$ & $\displaystyle\Big(\frac{1}{2},1 \Big)$
   %
   %
   &
   $ 0<q<1<p$  &
   \textit{super-lin / sub-lin}
   \\\\
   $\{1\}$ & $\{1\}$
  %
  %
  & $q=1<p=3$   &
  \textit{super-lin / lin}
   \\\\
  $(1,2)$ & $(1,+\infty)$
  %
  %
  &$1<q<p $  &
  \begin{minipage}{1.5 cm}
   \textit{\centerline{super-lin}}
   \end{minipage}
   \\\\
  $(2,+\infty)$ & $(-\infty,0)$
  %
  %
  &$p<q<0 $   & \textit{sing }
   \\\\
  \hline\hline
\end{tabular}
\end{table}
\end{small}

\begin{description}
  \item[$\underline{k = 2}$] so $0 < \displaystyle 2-k_{1}
= \frac{1}{2} < k_{1}=\frac{3}{2}$. It follows that \textsc{Table
\ref{tab:scalar curvature nonlinearity type m = 1, k = 2}}.
\end{description}

\begin{small}
\begin{table}
  \caption{\small Nonlinearities in scalar curvature equation type
\eqref{eq:scalar curvature m ge 2} for $m = 1 $ and $k = 2$, see
\textit{Notation \ref{nota: right columns in tables 45678}} }
\label{tab:scalar curvature nonlinearity type m = 1, k = 2}
\begin{tabular}{cccc}
\hline\hline
\\
  $\mu \in $ &$ \alpha \in$
   %
   %
   & $p,q$ &
    \textit{ type of $p,q$ non-linearity}
    %
   \\\\
 \hline
   \\
  $\displaystyle \Big(-\infty,-\frac{3}{2}\Big)$ &
$\displaystyle\Big(0,\frac{1}{3}\Big)$
   %
   %
   &$q<p<0$  & \textit{sing}
   \\\\
  $\displaystyle \Big\{-\frac{3}{2}\Big\}$ &
$\displaystyle\Big\{\frac{1}{3}\Big\}$
   %
   %
   &$\displaystyle q=-\frac{2}{3}<p=0$  &
   \textit{\underline{non}-hom / sing}
   \\\\
  $\displaystyle\Big(-\frac{3}{2},0\Big)$ &
$\displaystyle\Big(\frac{1}{3},\frac{2}{3}\Big)$
   %
   %
   &$q<0<p<1$  &
   \textit{sub-lin / sing}
   \\\\
$\displaystyle\Big(0,\frac{1}{2}\Big)$ &
$\displaystyle\Big(\frac{2}{3},1\Big) $
   %
   %
   &
   $ q<0<1<p$  &
   \textit{super-lin / sing}
   \\\\
   $\displaystyle\Big\{\frac{1}{2}\Big\}$ & $\{1\}$
   %
   %
   &
   $ q=0<p=2$  &
   \textit{super-lin / \underline{non}-hom}
   \\\\
   $\displaystyle\Big(\frac{1}{2},1\Big)$ & $(1,2)$
   %
   %
   &
   $ 0<q<1<p$  &
   \textit{super-lin / sub-lin}
   \\\\
   $\{1\}$ & $\{2\}$
   & $q=1<p=5$   &
   \textit{super-lin / lin}
   \\\\
   $\displaystyle\Big(1,\frac{3}{2}\Big)$ & $(2,+\infty)$
  %
  %
  &$1<q<p $  & \textit{super-lin }
   \\\\
   $\displaystyle\Big(\frac{3}{2},+\infty\Big)$ & $(-\infty,0) $
  %
  %
  &$p<q<0 $   & \textit{ sing }
   \\\\
\hline\hline

\end{tabular}
\end{table}
\end{small}

\bigskip
\begin{description}
  \item[$\underline{k = 1}$]in this case $0 < 2-k_{1}= k_{1}=1$. But
 since $S_{F}\equiv 0$, $q$ is non-influent. Thus we obtain \textsc{Table
\ref{tab:scalar curvature nonlinearity type m = 1, k = 1} }.
\end{description}

\begin{small}
\begin{table}
  \caption{\small Nonlinearities in scalar curvature equation type
\eqref{eq:scalar curvature m ge 2} for $m = 1 $ and $k = 1$, see
\textit{Notation \ref{nota: right columns in tables 45678}} }
%
\label{tab:scalar curvature nonlinearity type m = 1, k = 1}
\begin{tabular}{cccc}
\hline\hline
  \\
   $\mu \in $ &$ \alpha \in $
   & $p$ &
    \textit{ type of $p,q$ non-linearity}
    %
   \\\\
 \hline
   \\
 $(-\infty,-1)$ & $\displaystyle\Big(0,\frac{1}{2}\Big)$
   &$p<0$  & \textit{sing}
   \\\\
 $ \{-1\}$ & $\displaystyle\Big\{\frac{1}{2}\Big\}$
   %
   %
   &$p=0$  &
   \begin{minipage}{1.5 cm}
   \textit{\centerline{\underline{non}-hom}}
   \end{minipage}
   \\\\
   $(-1,0)$ & $\displaystyle\Big(\frac{1}{2},1\Big)$
   %
   %
   &$0<p<1$  & \textit{sub-lin}
   \\\\
   $(0,1)$ & $(1,+\infty)$
   %
   %
   &
   $ 1<p$  & \textit{ super-lin}
   \\\\
  $(1,+\infty)$ & $(-\infty,0) $
  %
  %
  &$p<0 $   & \textit{ sing }
   \\\\
  \hline\hline
\end{tabular}
\end{table}
\end{small}

\section{Some Examples and Final Remarks}
\label{sec:Some Examples} We consider the usual definition of
Einstein manifolds (see \cite{Aubin82,Aubin98,
Beem-Ehrlich-Easley96, Hawking-Ellis73,Lovelock-Rund75,ONeil83}).
For some other alternative but close definitions see
\cite{Besse87}. For dimension $\ge 3$ these definitions are
coincident.

\begin{defin}
A pseudo-Riemannian manifold $(N_{n},h)$ is said to be an Einstein
manifold with $\lambda \in C^\infty (N)$ if and only if
$Ric_{h}=\lambda h$.
\end{defin}

\noindent Thus, the followings hold by letting $(N_{n},h)$ be a
pseudo-Riemannian manifold,
\begin{description}
  \item[i] if $(N_{n},h)$ is Einstein with $\lambda$ and $n \ge
  3$, then $\lambda$ is constant and $\lambda = S_{N}/n$, where
  $S_{N}$ is the scalar curvature of $(N_{n},h)$.
  \item[ii] if $(N_{n},h)$ is Einstein with $\lambda$ and $n = 2$,
then $\lambda$ is not necessarily constant.
\end{description}


\begin{rem}\label{rem:Einstein:exam: base m=1,k=k}
Let $M=$ \twistpar{B_m}{F_k}{\psi^\mu}{\psi} be a
\bcwpar{\psi,\mu} such that the Ricci curvature tensor \ric ~ is
given by \eqref{eq:ricci special 1}.
%
%
So, M is an Einstein manifold with $\lambda$ if and only if
$(F,g_F)$ is Einstein with $\nu$ constant (note that when $k=2$,
$\nu$ is constant by the equations and not by the above item
\textbf{i}) and the system that follows is verified
\begin{equation}\label{eq:einstein: special m=m}
\begin{split}
&\displaystyle{\lambda \psi^{2\mu}g_B={\rm Ric}_B+
\beta^{H}\frac{1}{\psi^{\frac{1}{\alpha ^{H}}}}{\rm
H}_{B}^{\psi^{\frac{1}{\alpha ^{H}}}} - \beta^{\Delta}
\frac{1}{\psi ^{\frac{1}{\alpha ^{\Delta}}}} {\Delta}_{B}
\psi^{\frac{1}{\alpha ^{\Delta}}} }g_B
\textrm{ on } \mathcal{L}(B)\times \mathcal{L}(B)
%
%
\\
&\displaystyle{\lambda \psi^2=\nu
-\frac{1}{ \psi^{2(\mu - 1)}}
\frac{\beta^{\Delta}}{\mu} \frac{1}{\psi ^{\frac{1}{\alpha
^{\Delta}}}} {\Delta}_{B} \psi^{\frac{1}{\alpha ^{\Delta}}}
},
%
%
\end{split}
\end{equation}
where the coefficients are given by \eqref{eq:zeta eta HD
formulas}. Compare this system with the well known results for an
arbitrary warped product in \cite{Besse87, Kim-Kim01, ONeil83}.

\noindent So taking the $g_{B}-$trace of the first equation in
\eqref{eq:einstein: special m=m} results
\begin{equation}\label{eq:einstein: special trace m=m}
\begin{split}
&\displaystyle{\lambda m \psi^{2\mu}=S_B+
\beta^{H}\frac{1}{\psi^{\frac{1}{\alpha
^{H}}}}\Delta_{B}{\psi^{\frac{1}{\alpha ^{H}}}} - m\beta^{\Delta}
\frac{1}{\psi ^{\frac{1}{\alpha ^{\Delta}}}} {\Delta}_{B}
\psi^{\frac{1}{\alpha ^{\Delta}}} }
\\
&\displaystyle{\lambda\psi^{2}=\nu
-\frac{1}{ \psi^{2(\mu - 1)}}
\frac{\beta^{\Delta}}{\mu} \frac{1}{\psi ^{\frac{1}{\alpha
^{\Delta}}}} {\Delta}_{B} \psi^{\frac{1}{\alpha ^{\Delta}}}
}.
%
%
\end{split}
\end{equation}
At this point we observe that we meet all the hypothesis to apply
\textbf{Lemma \ref{lem:2.1}}, thus \eqref{eq:einstein: special
trace m=m} is equivalent to
\begin{equation}\label{eq:einstein: special trace m=m  2}
\begin{split}
&\displaystyle{\lambda m \psi^{2\mu}=S_B+
\beta_{tr}\frac{1}{\psi^{\frac{1}{\alpha
_{tr}}}}\Delta_{B}{\psi^{\frac{1}{\alpha _{tr}}}}  }
\\
&\displaystyle{\lambda\psi^{2}=\nu
-\frac{1}{ \psi^{2(\mu - 1)}}
\frac{\beta^{\Delta}}{\mu} \frac{1}{\psi ^{\frac{1}{\alpha
^{\Delta}}}} {\Delta}_{B} \psi^{\frac{1}{\alpha ^{\Delta}}}
},
%
%
\end{split}
\end{equation}
where
\begin{equation}\label{eq:eta beta trace}
  \begin{array}{l}
     \alpha_{tr} = \displaystyle\frac{\zeta_{tr}}{\eta_{tr}},
           \\
     \beta_{tr} = \displaystyle\frac{\zeta_{tr}^2}{\eta_{tr}},
  \end{array}
\end{equation}
with
\begin{equation}\label{eq:zeta eta trace}
  \begin{array}{lll}
     \zeta_{tr}
     &=& \displaystyle\frac{\beta^H}{\alpha^H}-m
\frac{\beta^\Delta}{\alpha^\Delta} =
     \zeta^H-m\zeta^\Delta
     =  -2(m-1)\mu-k ,
     \\
     \eta_{tr}
     &=& \displaystyle\frac{\beta^H}{(\alpha^H)^2}-m
\frac{\beta^\Delta}{(\alpha^\Delta)^2} =
     \eta^H-m\eta^\Delta
     \\
     &=&-(m-1)\mu[(m-2)\mu + k]+ k(\mu -1).
  \end{array}
\end{equation}

Note that for $m=1,$ we have $S_B \equiv 0$, thus the system
\eqref{eq:einstein: special trace m=m} and hence
\eqref{eq:einstein: special trace m=m  2} are equivalent to the
Einstein condition with $\lambda$. In this case, the coefficients
take the form
\begin{equation}\label{eq:alpha beta trace m=1}
  \begin{array}{lclcl}
     \alpha_{tr} &=& \displaystyle\frac{-1}{\mu - 1}, &&
            \\
     \beta_{tr} &=& \displaystyle\frac{k}{\mu - 1}, &&
            \\
     \alpha^{\Delta} &=& \displaystyle\frac{1}{-\mu + k}, &&
            \\
     \beta^{\Delta} &=& \displaystyle\frac{\mu}{-\mu + k}. &&
  \end{array}
\end{equation}
\end{rem}

\bigskip

\begin{exam}
\label{exam:+R1 times Einstein} First of all, note that the
interesting solutions of the involved ordinary differential
equations must be nonnegative and moreover positive for us. So along
this example, when we speak of solutions, it should be understood
that we consider only positive solutions, unless explicitly
mentioned otherwise. We now consider \textit{Remark
\ref{rem:Einstein:exam: base m=1,k=k}} with \emph{$B$ as a real
interval (i.e. $m =\dim B = 1$) equipped with the usual metric $\pm
dr^2$ and $(F^k,g_F)$ is an Einstein manifold with $\nu$}. We
immediately observe that in $(B,\pm dr^2)$, we have the following
expressions:

\begin{equation}\label{eq:metric dim 1 on R}
 \begin{split}
  &\nabla^B(\cdot)=\pm (\cdot)^{'},\\
  &|\nabla^B(\cdot)|_B^2=\pm  |(\cdot)^{'}|^2 ,\\
  &\displaystyle {\rm
H}_{B}^{(\cdot)}\Big(\frac{\partial}{\partial r
},\frac{\partial}{\partial r}\Big)=(\cdot)^{''},\\
  &\Delta_B(\cdot)= \pm (\cdot)^{''},
 \end{split}
\end{equation}

%
%

\noindent where $(\cdot)^{'}$ means the usual derivative with
respect to $r$. Thus, by \eqref{eq:einstein: special trace m=m 2},
if $\mu \in \mathbb{R}\setminus\{0,1,k,\overline{\mu}_\pm\}$
\footnote{The signs $\pm$ in $\mu_\pm$ are not relative to the
signs in the metric $\pm dr^2$, these are relative only with
\textsc{Table \ref{tab:ricci mu exceptional cases m 1,2}}.}
(see \textsc{Table \ref{tab:ricci mu exceptional cases m 1,2}})
the corresponding \bcwpar{\psi,\mu} is an Einstein manifold with
$\lambda$ if and only if $(\lambda,\psi)$ verifies the system

\begin{equation}\label{eq:einstein: exam: m=1,k=k}
\begin{split}
&\displaystyle{\lambda  \psi^{2\mu}= \pm
\frac{k}{\mu - 1}\frac{1}{\psi^{1-\mu}}(\psi^{1-\mu})^{''} }
\\
&\displaystyle{\lambda\psi^{2}=\nu
\mp \frac{1}{ \psi^{2(\mu - 1)}}
\frac{1}{-\mu + k} \frac{1}{\psi ^{-\mu + k}} (\psi^{-\mu +
k})^{''}
},
%
%
\end{split}
\end{equation}
or still by changing variables $v=\psi^{1-\mu}$, if and only if
$(\lambda,v)$ verifies the system
\begin{equation}\label{eq:einstein: exam: m=1,k=k v}
\begin{split}
(a)\hspace{1in}&\displaystyle{\lambda  v^{\frac{2\mu}{1-\mu}}=
\pm \frac{k}{\mu - 1}\frac{1}{v} v^{''}}
\\
(b)\hspace{1in}&\displaystyle{\lambda v^{\frac{2}{1-\mu}}=\nu
\mp v^2
\frac{1}{-\mu + k} \frac{ 1}{v ^{\frac{-\mu + k}{1-\mu}}}( v
^{\frac{-\mu + k}{1-\mu}})^{''}
}.
%
%
\end{split}
\end{equation}
So, applying \eqref{eq:lappotencia} to the right hand side of
\eqref{eq:einstein: exam: m=1,k=k v}$-(b)$ results that a solution
$(\lambda,v)$ of \eqref{eq:einstein: exam: m=1,k=k v}$-(a)$ is
solution of \eqref{eq:einstein: exam: m=1,k=k v}$-(b)$ if and only
if it is a solution to the first order ordinary differential
equation
\begin{equation}\label{eq:einstein: exam: m=1,k=k v:2}
\displaystyle{\lambda v^{\frac{2}{1-\mu}}=\nu
\mp
\frac{k-1}{(1-\mu)^2}  (v^{'})^2   + \frac{\lambda}{k}
v^{\frac{2}{1-\mu}}}
\end{equation}
or equivalently to
\begin{equation}\label{eq:einstein: exam: m=1,k=k v:3}
(k-1)\left(\pm \displaystyle{\frac{1}{(1-\mu)^2} (v^{'})^2 +
\frac{\lambda}{k} v^{\frac{2}{1-\mu}}}\right)=\nu.
%
\end{equation}

%
%
We divide the study in two cases, namely.
%
%

\begin{description}
  \item[\underline{$k \ge 2$}] In this case, Equation
\eqref{eq:einstein: exam: m=1,k=k v:3} is central, taking its
derivative we obtain that any regular solution of this verifies
\begin{equation}\label{eq:einstein: exam: m=1,k=k v-(a):2}
    \displaystyle 2\frac{k-1}{1-\mu}v ^{'}\left(\pm \frac{1}{1 - \mu }
v^{''} + \frac{\lambda}{k} v^{\frac{2}{1-\mu}-1}\right)=0 .
\end{equation}
Hence, we have the following result:

\medskip
\begin{minipage}{10 cm}
  \textit{If $(\lambda,v)$ is a solution of \eqref{eq:einstein: exam:
m=1,k=k  v},
  then it is a solution of  \eqref{eq:einstein: exam: m=1,k=k v:3}.
  Moreover, if $(\lambda,v)$
  is a solution of \eqref{eq:einstein: exam: m=1,k=k v:3} then $v$
  is constant or is a solution of \eqref{eq:einstein: exam: m=1,k=k
  v}.
  }
  \end{minipage}
\medskip


\noindent Thus, we have proved:

\medskip
\begin{minipage}{10 cm}
  \textit{If a \bcwpar{\psi,\mu} is Einstein with $\lambda,$
  then $0<v=\psi^{1-\mu}$ satisfies Equation
  \eqref{eq:einstein: exam: m=1,k=k v:3}, where $\lambda$ is
  necessarily
  constant (indeed $m+k\ge 3$). Furthermore, if
  $0<v=\psi^{1-\mu}$ is a nonconstant solution of Equation
  \eqref{eq:einstein: exam: m=1,k=k v:3}, then the corresponding
  \bcwpar{\psi,\mu} is Einstein with $\lambda$. Furthermore, if $0<\psi$ is a constant,
  then a
  \bcwpar{\psi,\mu} is Einstein  if and only if
  $\lambda=0=\nu$.}
\end{minipage}
\medskip

\noindent We observe that Equation \eqref{eq:einstein: exam:
m=1,k=k v:3} may be solved by the method of separation of
variables
\begin{equation}\label{eq:einstein: exam: m=1,k=k v:4}
\displaystyle{\frac{dv}{dr} = v^{'} =\sqrt{\pm (1-\mu)^2
\Big(\frac{\nu}{k-1} - \frac{\lambda}{k}
v^{\frac{2}{1-\mu}}\Big)}.
}
\end{equation}
Thus, its solutions are given by
\begin{equation}\label{eq:einstein: exam: m=1,k=k v:5}
\displaystyle{
\int^{v} \frac{1}{
\sqrt{\displaystyle \pm (1-\mu)^2 \Big(\frac{\nu}{k-1} -
\frac{\lambda}{k} w^{\frac{2}{1-\mu}}\Big)} }dw
= r .
}
\end{equation}
For suitable values of the parameters, the latter integral may be
solved by applying special functions (more specifically,
hypergeometric functions called also Gauss-Kummer series and
elliptic functions, see for example \cite{{BirkhoffRota78},
{Weinstein}} or apply Mathematica, Maple etc.). As we mentioned in
\S 1, metrics of this type are considered in Randal-Sundrum theory
\cite{Ito01} and in super-gravity theories.

One particular simpler case of the above results corresponds to
$\mu=-1$, namely.


\bigskip

\underline{\textbf{\bcwpar{\psi,-1} with $k \ge 2$}} In this case,
Equation \eqref{eq:einstein: exam: m=1,k=k v:4} reduces to
\begin{equation}\label{eq:einstein: exam: m=1,k=k,mu=-1 v:4}
\displaystyle{\frac{dv}{dr} = v^{'} =\sqrt{\pm 4
\Big(\frac{\nu}{k-1} - \frac{\lambda}{k} v\Big)}
}
\end{equation}
and \eqref{eq:einstein: exam: m=1,k=k v:5} to
\begin{equation}\label{eq:einstein: exam: m=1,k=k,mu=-1 v:5}
\displaystyle{
r = \int^{v} \frac{1}{
\sqrt{\displaystyle \pm 4 \Big(\frac{\nu}{k-1} - \frac{\lambda}{k}
w\Big)} }dw
= \mp \frac{k}{\lambda} \sqrt{\frac{v\lambda-kv\lambda+k\nu}{\mp k
\pm k^2}} + \gamma ,
}
\end{equation}
with a real constant $\gamma$. Hence,
\begin{equation}\label{eq:einstein: exam: m=1,k=k,mu=-1 v:6}
\displaystyle{
v(r) = \mp \frac{\lambda}{k} (r + \gamma)^2 +  \frac{\nu}{\lambda}
\frac{k}{k-1}.
}
\end{equation}

\end{description}

\medskip

\begin{description}
  \item[\underline{$k = 1$}] First of all, $\nu = 0$ and $\lambda \in
C^\infty(B)$. Hence, unlike to the case of $k \ge 2$,
\eqref{eq:einstein: exam: m=1,k=k v:3} gives no information and
 \eqref{eq:einstein: exam: m=1,k=k
  v}$-(a)$ and \eqref{eq:einstein: exam: m=1,k=k v}$-(b)$ coincide.
Thus, we proved that:

\medskip
\begin{minipage}{10 cm}
  \textit{A \bcwpar{\psi,\mu} is Einstein with $\lambda \in
C^\infty(B)$ if and only if $0<v=\psi^{1-\mu}$ satisfies
\eqref{eq:einstein: exam: m=1,k=k
  v}$-(a)$ with $k=1$}.
  \end{minipage}

\bigskip

For the completeness of the exposition, we will write a few lines
about possibly the easiest case that follows.

\medskip

 \underline{\textbf{\bcwpar{\psi,-1} with $k=1$}} Here,
  \eqref{eq:einstein: exam: m=1,k=k  v}$-(a)$ takes the
 trivial form
\begin{equation}\label{eq:einstein: exam: m=1,k=1 v:1}
v^{''} = \mp 2\lambda,
\end{equation}
where $\lambda \in C^\infty(B)$. So
\begin{equation}\label{eq:ricci special eins 1.4}
v(r) = \mp 2\int^{r} \int ^{\omega}\lambda(\tau)d\tau d\omega.
\end{equation}

In particular, if $\lambda$ is constant then this results
$\psi^{2} (r) = v(r) = \mp 2\lambda r^{2} + a r + b$, with $a$ and
$b$ real constants such that $\mp 2\lambda r^{2} + a r + b$ is
positive. It is clear that the latter condition depends on the
base interval $B_{1}$ and the parameter $\lambda$.
\end{description}
\end{exam}


\begin{exam}
\label{exam:+compact 1 times Einstein} Like in \textit{Example
\ref{exam:+R1 times Einstein}}, we consider only positive
solutions, unless otherwise explicitly mentioned. By applying
\textit{Remark \ref{rem:Einstein:exam: base m=1,k=k}} when $B$ is
a compact Riemannian manifold of $\dim B = m = 1$ with metric
$g_B$ and $(F^k,g_F)$ is an Einstein manifold with $\nu$, we have
that if $\mu \in \mathbb{R}\setminus\{0,1,k,\overline{\mu}_\pm\}$
(see \textsc{Table \ref{tab:ricci mu exceptional cases m 1,2}})
the corresponding \bcwpar{\psi,k} is an Einstein manifold with
$\lambda$ if and only if $(\lambda,\psi)$ verifies the system
\begin{equation}\label{eq:einstein: +compact 1 times Einstein}
\begin{split}
&\displaystyle{\lambda  \psi^{2\mu}=
\frac{k}{\mu - 1}\frac{1}{\psi^{1-\mu}}{\Delta_B(\psi^{1-\mu})} }
\\
&\displaystyle{\lambda\psi^{2}=\nu
-\frac{1}{ \psi^{2(\mu - 1)}}
\frac{1}{-\mu + k} \frac{1}{\psi ^{-\mu + k}} \Delta_B(\psi^{-\mu
+ k})
},
%
%
\end{split}
\end{equation}
Thus, by integrating on $B$ and applying the compactness of $B$
and also considering the positivity of $\psi$ we conclude that
$\lambda = \nu = 0$ and $\psi$ is a positive constant. So, we
proved that:

\medskip
\begin{center}
\begin{minipage}{10 cm}
\textit{Let $B$ be a compact Riemannian manifold of $\dim B = m =
1$ with metric $g_B$ and $(F^k,g_F)$ be an Einstein manifold with
$\nu$ where $\mu \in
\mathbb{R}\setminus\{0,1,k,\overline{\mu}_\pm\}$. A
\bcwpar{\psi,\mu} is Einstein with $\lambda$ if and only if
$\lambda = \nu = 0$ and $\psi$ is a positive constant (in
particular, a trivial product).}
\end{minipage}
\end{center}
\medskip
\end{exam}

\begin{rem} \label{rem: Einstein base flat}The same order of ideas of
\textit{Example \ref{exam:+compact 1 times Einstein}} and
considering especially \eqref{eq:einstein: special trace m=m 2}
and \eqref{eq:einstein: special m=m}, allow us to prove the
following:
\medskip
\begin{center}
\begin{minipage}{10 cm}
\textit{Let $(B_m,g_B)$ be a scalar flat compact Riemannian
manifold and $(F_k,g_F)$ be a pseudo-Riemannian manifold.
Furthermore, suppose that $\mu \in \mathbb{R}
\setminus\{0,1,k,\overline{\mu}_\pm\}$. A \bcwpar{\psi,\mu} is
Einstein with a constant $\lambda$ if and only if $(F_k,g_F)$ is
Einstein with $\nu = 0$, $\lambda=0$ and $\psi$ is a positive
constant (in particular a usual product) and $(B_m,g_B)$ is
Ricci-flat.}
\end{minipage}
\end{center}

\end{rem}

\begin{rem}\label{rem: Einstein functional condition} Let $k \ge
2$ be and let us assume
the hypothesis of Remark \ref{rem:Einstein:exam: base m=1,k=k}.
\begin{description}
  \item[i]
It is easy to verify that in Equation \eqref{eq:einstein: special
trace m=m 2},
\begin{equation}\label{eq:m=2;H-JS;0}
  \alpha_{tr}=\alpha^{\Delta}
\end{equation}
if and only if
\begin{equation}\label{eq:m=2;H-JS;1}
  (m-1)(m-2)\mu^2 + 2(m-1)k\mu + k(k-1)=0.
\end{equation}
Note that the latter equation \eqref{eq:m=2;H-JS;1} is also
equivalent to
\begin{equation}\label{eq:m=2;Delta=H}
 \alpha^{\Delta}=\alpha^{H}.
\end{equation}
 Since $m,k \in \mathbb{N}$ and $k \ge 1$, for any
$m
> 2$ the equation \eqref{eq:m=2;H-JS;1} has two real solutions,
namely
\begin{equation}\label{eq:m=2;H-JS;2}
 \tilde{\mu}_{\pm}=\frac{-(m-1)k \pm \sqrt{
 (m-1)k(k+m-2)
 }}{(m-1)(m-2)
 },
\end{equation}
while for $m=2$ has only one solution
\begin{equation}\label{eq:m=2;H-JS;3}
 \tilde{\mu}=\frac{1-k}{2}.
\end{equation}
We remark here that the latter is exactly the value of the
parameter considered by H-J. Schmidt in his studies about
Birkhoff's theorems in \cite{Schmidt97} (see \textbf{vi} in \S 1).

  \item[ii]
  If Equation \eqref{eq:m=2;H-JS;0} is satisfied for some $\mu$,
  then \eqref{eq:einstein: special trace m=m  2} implies the functional
equation
\begin{equation}\label{eq:m=2;H-JS;4}
  \lambda m \psi^{2\mu}=S_B + \beta_{tr} (\nu - \lambda
  \psi^2)\psi^{2(\mu-1)}\frac{\mu}{\beta^\Delta},
\end{equation}
or equivalently, by \eqref{eq:zeta eta HD formulas},
\eqref{eq:zeta eta HD}, \eqref{eq:eta beta trace} and
\eqref{eq:zeta eta trace}
\begin{equation}\label{eq:m=2;H-JS;5}
  \lambda m \psi^{2\mu}=S_B - [2(m-1)\mu + k] (\nu - \lambda
  \psi^2)\psi^{2(\mu-1)},
\end{equation}
or still
\begin{equation}\label{eq:m=2;H-JS;6}
  [m-2(m-1)\mu - k]\lambda  \psi^{2\mu} + [2(m-1)\mu + k] \nu
\psi^{2(\mu-1)}=S_B.
\end{equation}

We observe that if $\mu$ is such that \eqref{eq:m=2;H-JS;0} is
satisfied, we reobtained \textit{Remark \ref{rem: Einstein base
flat}} (for this specific value of $\mu$) \textit{without} the
hypothesis of compactness of the base, as a consequence of
\eqref{eq:m=2;H-JS;6} and \eqref{eq:einstein: special m=m}.

  \item[iii] When $m=2$ and $\mu $ is like in
  \eqref{eq:m=2;H-JS;3}, then \eqref{eq:m=2;H-JS;6} takes the form
\begin{equation}\label{eq:m=2;H-JS;7}
  \lambda  \psi^{1-k} +  \nu \psi^{-(k+1)}=S_B.
\end{equation}

\end{description}
\end{rem}

\begin{exam} Now, we consider an interesting application of
\eqref{eq:einstein: special trace m=m  2} with $m=2$ and $k \ge
2$, containing as particular case the Schwarzchild type metrics
considered in \textbf{i} of \S 1. Along the development of this
example, we will prove the statement that follows:
\medskip
\begin{center}
\begin{minipage}{10 cm}
  \textit{Let $(F_k,g_F)$ be Einstein with constant
  Ricci
  curvature
  $\nu$ and dimension $k \ge 2$. Then,
  $\mathbb{R}_{+} \times \mathbb{R} \times F_k$ furnished with a
metric
\begin{equation}\label{eq:H-J S gen}
\displaystyle{
  g=s^{\frac{1}{k}-1}\Big[ \frac{1}{4\sqrt{s}u^2(\sqrt{s})} ds^2 \pm
4 \sqrt{s} u^2(\sqrt{s}) dy^2  \Big]
  + s^{\frac{2}{k}} g_F
  },
\end{equation}
is Einstein with constant
Ricci curvature $\lambda$ where $(s,y)\in \mathbb{R}_{+} \times
\mathbb{R}$ if and only if $u^2$ is given by \eqref{eq:einstein:
special m=2 k, 6} below, where $\lambda $ and $C$ are such that the
right hand side of \eqref{eq:einstein: special m=2 k, 6} results
positive. }
\end{minipage}
\end{center}
\medskip

Let $(F_k,g_F)$ be Einstein with $\nu$ and
$(B_2,g_B)=(\mathbb{R}_{+} \times \mathbb{R}, g_B)$ be a
pseudo-Riemannian manifold endowed with the metric
\begin{equation}\label{eq:m=2=k;H-JS;1}
  g_B=(\psi_1(s))^{2(-1)}ds^2 \pm (\psi_1(s))^{2}dy^2,
\end{equation}
where $\psi_1 $ is defined as $\psi_1(s)=2 s^\frac{1}{4}
u(s^\frac{1}{2})$, like in \eqref{eq:Schwarzschild to bcwp-3}. So by
applying the second row of \textsc{Table \ref{tab:scalar curvature
nonlinearity type m = 1, k = 1}} we have,
\begin{equation}\label{eq:m=2=k;H-JS;2}
  S_B(s)=-\Delta_{ds^2}\psi_1^{2}(s)=\mathcal{S}_B
  u^{2}|_{r=s^{\frac{1}{2}}},
\end{equation}
where $\mathcal{S}_B$ is the linear second order ordinary
differential operator defined by
\begin{equation}\label{eq:m=2=k;H-JS;3}
  \mathcal{S}_B f (r)= r^{-3}f(r) - r^{-2}\frac{d}{dr}f\Big|_{r} -
  r^{-1}\frac{d^2}{dr^2}f\Big|_{r}, f \in C^\infty(\mathbb{R}_+).
\end{equation}

We now consider $B_2 \times F_k = \mathbb{R}_{+} \times \mathbb{R}
\times F_k$ endowed with the metric
\begin{equation}\label{eq:m=2 k;H-JS;4}
g =(\psi_2(s,y))^{2 \mu_2}g_B + (\psi_2(s,y))^{2}g_F,
\end{equation}
where $\psi_2(s,y)=s^{\frac{1}{k}}$ and
$\displaystyle\mu_2=\frac{1-k}{2}$ (compare with
\eqref{eq:Schwarzschild to bcwp-4} when $k=2$).

\noindent
Hence, since $(B_2 \times F_k,g)$ satisfies the hypothesis of
Remark \ref{rem: Einstein functional condition} (see Remark
\ref{thm:ricci sbcwp m = 1, 2}), if $(B_2 \times F_k,g)$ is
Einstein with $\lambda$, then $\psi_2$ satisfies
\eqref{eq:einstein: special trace m=m 2}. Furthermore, the
relation \eqref{eq:m=2;H-JS;0} is verified with $\displaystyle
\mu_2 = \tilde{\mu}$ (see \eqref{eq:m=2;H-JS;3}). Consequently,
$\psi_2$ must verify \eqref{eq:m=2;H-JS;7}, precisely
\begin{equation}\label{eq:m=2 k;H-JS;5}
  \lambda  \psi_2^{1-k} +  \nu \psi_2^{-(k+1)}=S_B.
\end{equation}
Therefore, by \eqref{eq:m=2=k;H-JS;2} and the definition of
$\psi_2$
\begin{equation}\label{eq:m=2 k;H-JS;6}
  \lambda  r^{\frac{2}{k}-2} +  \nu r^{-2- \frac{2}{k}}=\mathcal{S}_B
  u^{2}|_{r},
\end{equation}
or equivalently
\begin{equation}\label{eq:m=2 k;H-JS;7}
  \lambda  r^{1+\frac{2}{k}} +  \nu r^{1-\frac{2}{k}}=r^3
\mathcal{S}_B
  u^{2}|_{r},
\end{equation}
where $r=s^{\frac{1}{2}}$.

\noindent Note that the latter is an Euler (also called
equidimensional) equation. It is easy to show that for any real
constants $\nu$ and $\lambda$, the general solution of
\eqref{eq:m=2 k;H-JS;7} has the form
\begin{equation}\label{eq:m=2 k;H-JS;7 bis}
  u^2(r)= \lambda \Big(1-\Big(1+\frac{2}{k}\Big)^{2}\Big)^{-1}
r^{1+\frac{2}{k}}
  + \nu \Big(1-\Big(1-\frac{2}{k}\Big)^{2}\Big)^{-1}
r^{1-\frac{2}{k}} + v_h(r),
\end{equation}
where $v_h$ is a solution of the homogeneous equation
\begin{equation}\label{eq:m=2 k;H-JS;8}
  0 =u^2 - r \frac{d}{dr} u^2 \Big|_{r} -
  r^{2} \frac{d^2}{dr^2} u^2 \Big|_{r},
\end{equation}
namely a linear combination of $r$ and $r^{-1}$.

\noindent It is clear that the choices of $\nu, \lambda$ and $v_h$
will be such that the function $u^2$ be nonnegative.

\noindent Furthermore, we observe that among all the solutions of
\eqref{eq:m=2 k;H-JS;6} there are spurious solutions of
\eqref{eq:einstein: special m=m}, the reason is that
\eqref{eq:einstein: special trace m=m} is only a necessary condition
of \eqref{eq:einstein: special m=m}. Indeed, \eqref{eq:m=2 k;H-JS;7
bis} is a solution of \eqref{eq:einstein: special m=m} if and only
if $\displaystyle v_h(r)= C \frac{1}{r}$, where $C$ is an arbitrary
constant. 
In order to prove this, we note the following facts about
$(B_2,g_B)$ which is assumed as above:
\begin{description}
  \item[i] Since $m=2$,
\begin{equation}\label{eq:s-ric}
  \ric_B = \displaystyle \frac{1}{2} S_B g_B
\end{equation}
  \item[ii]By \emph{Proposition} \ref{prp: hessian bcwp},
\begin{equation}\label{eq:s-Hess}
  {\rm H}_{B}^s = \Big(\displaystyle \frac{1}{2} \Delta_{g_{B}} s\Big)
  g_B.
\end{equation}
  \item[iii] By \emph{Proposition} \ref{prp: laplacian bcwp} and
  the definition of $\psi_1$,
\begin{equation}\label{eq:s-Laplacian}
  \Delta_{g_{B}} s = 2 \psi_1 \frac{d}{ds}\psi_1 =
  \frac{d}{ds} \psi_1^2 = \mathcal{L}\Big(r,\displaystyle
\frac{d}{dr}\Big)u^2 \Big|_{s^{\frac{1}{2}}},
\end{equation}
where
\begin{equation}\label{eq:s-Laplacian-2}
  \mathcal{L}\Big(r,\displaystyle \frac{d}{dr}\Big)f \Big|_r =
2 \Big[r^{-1} + \frac{d}{dr} \Big] f \Big|_r.
\end{equation}

\end{description}

For $(B_2 \times F_k,g)$, since the coefficients given by
\eqref{eq:zeta eta HD formulas} verify \eqref{eq:m=2;H-JS;0}, they
take the values:
\begin{equation}\label{eq:alphaDH betaDH m=2 k}
  \begin{array}{l}
    \alpha^{\Delta} =   \alpha^{H} = \displaystyle \frac{1}{k}      \\
    \beta^{\Delta} = \displaystyle \frac{\tilde{\mu}}{k}   =
\frac{1-k}{2k}   \\
    \beta^{H} =  \displaystyle \frac{k}{2\tilde{\mu} - 1} = -1 ,
  \end{array}
\end{equation}
hence, \eqref{eq:einstein: special m=m} takes the form
\begin{equation}\label{eq:einstein: special m=2 k}
\begin{split}
&\displaystyle{\lambda \psi_2^{1-k}g_B={\rm Ric}_B
-\frac{1}{\psi_2^{k}}\Big({\rm H}_{B}^{\psi_2^{k}} +
\frac{1-k}{2k} {\Delta}_{B} \psi_2^{k} }g_B\Big)
\textrm{ on } \mathcal{L}(B)\times \mathcal{L}(B)
\\
&\displaystyle{\lambda \psi_2^2=\nu
-\psi_2^{k+1}
\frac{1}{k} \frac{1}{\psi_2 ^{k}} {\Delta}_{B} \psi_2^{k}
}.
\end{split}
\end{equation}
So by the definition of $\psi_2$, \eqref{eq:s-ric} and
\eqref{eq:s-Hess}, \eqref{eq:einstein: special m=2 k} results
equivalent to
\begin{equation}\label{eq:einstein: special m=2 k, 2}
\begin{split}
&\displaystyle{\lambda \psi_2^{1-k} = \frac{1}{2} S_B
-\frac{1}{2} \frac{1}{k} \frac{1}{\psi_2^k}  \Delta_{B} \psi_2^k
}
\\
&\displaystyle{\lambda \psi_2^{1-k}-\nu \psi_2^{-(k+1)}
=-
\frac{1}{k} \frac{1}{\psi_2 ^{k}} {\Delta}_{B} \psi_2^{k}
},
\end{split}
\end{equation}
or moreover, by easy computations, to
\begin{equation}\label{eq:einstein: special m=2 k, 3}
\begin{split}
&\lambda \psi_2^{1-k} + \nu \psi_2^{-(k+1)}= S_B
\\
&\lambda \psi_2 - \nu \psi_2^{-1} =
- \frac{1}{k}  {\Delta}_{B} \psi_2^{k} .
\end{split}
\end{equation}

Note in the above steps the reduction from $4$ to $2$ equations.
Furthermore the first equation of \eqref{eq:einstein: special m=2
k, 3} is exactly \eqref{eq:m=2 k;H-JS;5}. Recalling again that
$\psi_2(s,y)=s^{\frac{1}{k}}$, \eqref{eq:einstein: special m=2 k,
3} takes the form
\begin{equation}\label{eq:einstein: special m=2 k, 4}
\begin{split}
&\lambda s^{\frac{1}{k}-1} + \nu s^{-\frac{1}{k}-1}= S_B
\\
&\lambda s^{\frac{1}{k}} - \nu s^{-\frac{1}{k}} =
- \frac{1}{k}  {\Delta}_{B} s  ,
\end{split}
\end{equation}
and since $s^{\frac{1}{2}}=r$,
by \eqref{eq:m=2=k;H-JS;2} and \eqref{eq:s-Laplacian},
\begin{equation}\label{eq:einstein: special m=2 k, 5}
\begin{split}
&(a) \qquad \lambda r^{\frac{2}{k}-2} + \nu r^{-\frac{2}{k}-2}=
\mathcal{S}_B u^{2}|_{r}
\\
&(b) \qquad \lambda r^{\frac{2}{k}} - \nu r^{-\frac{2}{k}} =
- \frac{1}{k}  \mathcal{L}\Big(r,\displaystyle
\frac{d}{dr}\Big)u^2 \Big|_{r}.
\end{split}
\end{equation}
We observe that deriving the second equation of \eqref{eq:einstein:
special m=2 k, 5} and multiplying by $r^{-1}$, we obtain the first
equation. So any regular solution of $(b)$ is a solution of $(a)$ in
\eqref{eq:einstein: special m=2 k, 5}.

On the other hand it is easy to show that a general solution of
\eqref{eq:einstein: special m=2 k, 5}-$(b)$ is
\begin{equation}\label{eq:einstein: special m=2 k, 6}
  u^2(r)= \lambda \Big(1-\Big(1+\frac{2}{k}\Big)^{2}\Big)^{-1}
r^{1+\frac{2}{k}}
  + \nu \Big(1-\Big(1-\frac{2}{k}\Big)^{2}\Big)^{-1}
r^{1-\frac{2}{k}} + C \displaystyle
  \frac{1}{r},
\end{equation}
where $C$ is an arbitrary constant.
%

Thus, since \eqref{eq:einstein: special m=2 k, 5}-$(a)$ coincides
with \eqref{eq:m=2 k;H-JS;6}, a solution \eqref{eq:m=2 k;H-JS;7
bis} of the latter is a solution of \eqref{eq:einstein: special
m=2 k, 5} if and only if $\displaystyle v_h(r)= C \frac{1}{r}$,
where $C$ is an arbitrary constant Q.E.D.

\bigskip

As we mentioned in the first paragraph of this example, important
solutions of the Einstein vacuum equations are included in the above
discussion (namely, compare with \S1 \textbf{i}). We will write
explicitly some cases with $k=2$ but the situation is more general.
Let $B_2 \times F_2$ be endowed with a metric of the form
\begin{equation}\label{eq:m=2 k=2;H-JS;14}
g =s^{-\frac{1}{2}}g_B + s g_{F_2},
\end{equation}
where $(B_2,g_B)$ is like in \eqref{eq:m=2=k;H-JS;1}, $(s,y)\in
B_2=\mathbb{R}_{+} \times \mathbb{R}$ and $(F_2,g_F)$ is a
pseudo-Riemannian manifold of dimension $k=2$.
\begin{description}
  \item[Ricci flat]
  If $\lambda = 0$, then \eqref{eq:m=2 k;H-JS;7}
  takes the form
  \begin{equation}\label{eq:m=2=k;H-JS;8}
  \nu = u^2 - r(u^2)^{'} -
  r^{2}(u^2)^{''}.
\end{equation}
It is easy to verify that $u^2(r)= \nu + C \displaystyle
\frac{1}{r}$ is a solution of \eqref{eq:m=2=k;H-JS;8}. In
particular, when $C=-2M$, $M>0$ and $\nu=1$ we obtain the classical
Schwarzchild solution (compare with \eqref{eq:Schwarzschild-3} and
\eqref{eq:einstein: special m=2 k, 6}). While, the condition ``$C=0$
and $\nu=1$" arises the Minkowski metric of an empty space-time in
spherical terms.

  \item[Riemman-Schwarzchild] If $\lambda = -3$ and $\nu = 1$, then
\eqref{eq:m=2 k;H-JS;7}
  takes the form
  \begin{equation}\label{eq:m=2=k;H-JS;10}
  -3 r^2 + 1 = u^2 - r(u^2)^{'} -
  r^{2}(u^2)^{''}.
\end{equation}
It is easy to verify that for any positive $M$, $\displaystyle
u^2(r)=1-\frac{2M}{r}+r^2$ is a solution of \eqref{eq:m=2=k;H-JS;10}
(compare with \eqref{eq:Schwarzschild-1b} and \eqref{eq:einstein:
special m=2 k, 6}).
\end{description}

\noindent Thus, Equation \eqref{eq:m=2 k;H-JS;7} contains a large
family of important solutions of the Einstein equation. An
analogous procedure can be applied to build the static BTZ (2+1)-
black hole solution, we leave the computations to the reader (see
\cite{Aharony-Gubser-Maldacena-Ooguri-Oz,
Banados-Teitelboim-Zanelli92,
Banados-Henneaux-Teitelboim-Zanelli93,
DobarroUnal04,Hong-Choi-Park03} for details about BTZ).


\end{exam}

\bigskip

\begin{rem}
Let $F = (F_k,g_{F})$ be a pseudo-Rieman\-nian Einstein manifold
with constant $\nu$ and dimension $k \geq 1$.

\noindent We recall the principal result in \cite{Kim-Kim01} in
the context of Riemannian manifolds, namely: an Einstein warped
product with a non-positive scalar curvature and \textit{compact}
base is a trivial Riemannian product space, so that the warping
function results constant. Thus, if $B = (B_m,g_{B})$ is a
\textit{compact} Rieman\-nian manifold with dimension $m \geq 3$
and $\displaystyle \mu \in \mathbb{R}\setminus \{0, 1,
\overline{\mu}, \overline{\mu}_{\pm}\}$ (compare with Theorem
\ref{thm:ricci sbcwp m ge 3}), then our system \eqref{eq:einstein:
special m=m} admits a non-constant positive solution only if
$\lambda
> 0$.
But if we let $F$ and $\mu$ be as above, then there exists a
metric on $B_m$ admitting no $\psi \in C^\infty_{>0}(B)$ such that
the corresponding \bcwpar{\psi,\mu} is Einstein with $\lambda >0$.
Indeed, multiplying the first equation of \eqref{eq:einstein:
special trace m=m  2} by $\displaystyle \psi^{\frac{1}{\alpha
_{tr}}} $ and integrating on $B$ respect to the measure $dg_B$
there results
\begin{equation*}
\lambda m \int_B \psi^{2\mu + \frac{1}{\alpha_{tr}}} dg_B = \int_B
S_B \psi^{\frac{1}{\alpha _{tr}}}dg_B.
\end{equation*}
Now, we recall the Aubin result ``any manifold of dimension $\ge
3$ possesses a complete metric of constant negative scalar
curvature" (see \cite{Aubin70,BlandKalka89,Lohkamp94}). So if
$g_B$ is a such metric on our compact $B_m$, i.e. $S_B < 0$, then
$\lambda $ cannot be positive (contradiction).

In conclusion, let $F$ and $\mu$ as above. On every compact
manifold $B$ of dimension $\ge 3$, there exits a Riemannian metric
$g_B$ such that a \bcwpar{\psi,\mu} with base $(B,g_B)$ is
Einstein with $\lambda $ if and only if $\psi$ is constant,
$(B,g_B)$ is Einstein with $\lambda m \psi^{2 \mu}$ and $\lambda
\psi^2  = \nu \le 0$.


The case $\mu = 0$, i.e. singly warped product, was considered in
\cite{Mustafa05}. The remaining values of $\mu$ (i.e. $1,
\overline{\mu}, \overline{\mu}_{\pm}$) can be analyzed with an
analogous approach with suitable changes, yet by applying
\eqref{eq:ricci special 4} and \eqref{eq:lappotencia}.


A particular example of the latter results (i.e. $\mu = -1$) is
the following interesting application of them:
%
%
%
Let $(F_k,g_{F})$ be a pseudo-Rieman\-nian Einstein manifold of
dimension $k \geq 1$. Then on any compact manifold $B_m$ of
dimension $\ge 3$ there exists a metric $g_B$ such that there is
no $\psi \in C^{\infty}_{>0}(B)$ such that $(B \times F,
\psi^{-2}g_B + \psi^2 g_F)$ is a non trivial (i.e $\psi $ non
constant) Einstein manifold.


\end{rem}

\section{Conclusions and future directions}
\label{sec:Conclusions}

Now, we would like to summarize the content of the paper and to
propose our future plans on this topic.

In brief, we introduced and studied curvature properties of a type
of product of two pseudo-Riemannian manifolds called \textit{base
conformal warped product} by us, roughly speaking the metric of a
such product is a mixture of a conformal metric on the base and a
warped metric.
%
%
%
%
As we mentioned in \S 1, these kind of metrics and considerations
about their curvatures are very frequent in different physical
areas, for instance relativity, extra-dimension theories
(Kaluza-Klein, Randall-Sundrum), string and super-gravity
theories; also in global analysis for example in the study of the
spectrum of Laplace-Beltrami operators on $p$-forms, etc.

In \S 2, we started our discussion by considering particular
families of either scalar or tensorial nonlinear partial
differential operators on pseudo-Riemannian manifolds and studied
useful identities verified by them. The latter allowed us to find
reduced expressions of the Ricci tensor and scalar curvature used
not only in \S 4 and \S 5, but also in the study of multiply
warped products in \cite{DobarroUnal04}. The operated reductions
can be considered as generalizations of those used by Yamabe in
\cite{Yamabe63}, in order to obtain the famous expression
\eqref{eq:yamabe} for the behavior of the scalar curvature under a
conformal change and those used in \cite{DobarroLamiDozo87} with
the same aim but for a singly warped product (see also Remark
\ref{rem:Lelong-Ferrand} for other particular application).

In \S 3, we defined precisely \textit{base conformal warped
products} of pseudo-Riemannian manifolds and computed their
Levi-Civita connection, Hessian, Laplace-Beltrami operator and
Riemannian curvatures.

In \S 4 and from then on, we concentrated on a very commonly used
physical ansatz, namely when the conformal factor acting on the
metric of the base and the warping function acting on the metric
of the fiber are related by an exponent, so that one is a power of
the other (see the examples in \S 1). We called a product manifold
furnished with a metric form like above as a special base
conformal warped product. Then, we turned our attention to the
structure of the relations that connect the different types of
curvatures, especially Ricci and scalar. More explicitly, we
obtained more approachable relations by applying the results of \S
2 but also some formulas even in some exceptional cases
corresponding to the situations where the results of \S 2 are
unapplicable.

In \S 5, we focused on a classification of the type of
nonlinearities arose in the relation among the involved scalar
curvatures of a \textit{special base conformal warped product},
previously obtained in \S 4. Similar to the study made in the
latter, we classified the nonlinearities according to the value of
%
%
the exponent parameter $\mu$, the dimensions of the base and the
fiber and finally
the scalar curvature of the fiber. The aim of this classification
is to study in future works the problem of prescribing
constant/nonconstant scalar curvature in \textit{special base
conformal warped products}, indeed in these problems, the type of
nonlinearities, ellipticity/hyperbolicity of the linear part of
differential equations connecting the involved scalar curvatures
and compactness of the base play a very central role.

At this point, we would like to note that the previous problems as
well as the study of the Einstein equation on \textit{base
conformal warped products}, \textit{special base conformal warped
products} and their \textit{generalizations to multi-fiber cases},
give rise to a reach family of interesting problems not only in
differential geometry and physics (see for instance, the several
recent works of R. Argurio, J. P. Gauntlett, S. Kachru, M. O.
Katanaev, J. Maldacena, H. -J. Schmidt, E. Silverstien, A.
Strominger, P. S. Wesson among many others), but also in non
linear analysis (see the different works of A. Ambrosetti, T.
Aubin, Y. Choquet-Bruat, J. F. Escobar, E. Hebey, R. Schoen, S.
-T. Yau among others), which will be the subject matter of future
works (see \cite{DobarroUnal07}).

In \S 6, we analyzed, investigated and characterized possible
solutions for the conformal and warping factors of a
\textit{special base conformal warped product} which guarantee
that the corresponding product
%
%
is Einstein. We apply the same order of ideas to a generalization
of the Schwarzchild metrics also. Among the considered cases there
are important metrics in questions of relativity, cosmology, hight
energy physics, etc.

\bigskip

\appendix{}
\numberwithin{equation}{section}
\renewcommand{\theequation}{\Alph{section}.\arabic{equation}}
\section{ }

\label{sec:appendix}

We first show some interesting properties about the behavior of
the Laplace\\-Beltrami operator under a conformal change in the
metric.

Let $N = (N_n,h)$ be a pseudo-Riemannian manifold of dimension $n$
and let
\begin{equation}
\Delta_{h}(\cdot) = \frac{1}{\sqrt{h}}
\partial_{i}(\sqrt{|h|}h^{ij}\partial_{j}(\cdot)),
\end{equation} be the Laplace-Beltrami operator related to the metric
$h$, where we denote the usual volume element by$\sqrt{|h|} :=
\sqrt{|det \, h|}$.

\begin{lem}
 \label{lem:conformal-laplacian-bis}
 Let $u \in C^{\infty}_{>0}(N)$ and $r \in \mathbb{R}$. Then,
 \begin{equation}
  u^{r}\Delta_{u^{r}h}(\cdot) =
  r {\frac{n-2}{2}} \, h\left(\frac{\nabla u}{u},\nabla (\cdot)
\right)
  + \Delta_{h}(\cdot).
  \label{eq:lap-conf4-ge3-bis}
\end{equation}
\end{lem}

\begin{proof} Denote $\tilde{h}=u^{r}h$, there results
$\tilde{h}_{ij} = u^{r} h_{ij}$, $\tilde{h}^{ij} = u^{-r} h^{ij}$
and  $\det \, \tilde{h}= u^{nr} \det \,h $. Thus,

\begin{equation}
  \begin{array}{lll}
    \displaystyle \Delta_{\tilde{h}}(\cdot) &=& \displaystyle
\frac{1}{u^{\frac{n}{2} r}\sqrt{h}}
\partial_{i}( u^{\frac{n}{2} r} \sqrt{h} \, u^{-r}
h^{ij}\partial_{j}(\cdot)) \\
    &=& \displaystyle \frac{1}{u^{\frac{n}{2} r}\sqrt{h}}
    \left[\left({\frac{n}{2}-1}\right)r u^{(\frac{n}{2}-1) r - 1} \partial_{i} u
\sqrt{h} \, h^{ij}\partial_{j}(\cdot) \right.\\
&&\left.+
u^{(\frac{n}{2}-1) r }
\partial_{i}(\sqrt{h} \, h^{ij} \partial_{j}(\cdot))\right].
  \end{array}
\end{equation}
So multiplying by $u^{r}$,
\begin{equation}
  u^{r}\Delta_{\tilde{h}}(\cdot) =
  \left({\frac{n}{2}-1}\right)r u^{-1} \partial_{i} u
h^{ij}\partial_{j}(\cdot)
  + \Delta_{h}(\cdot).
  \label{eq:lap-conf4-ge3-bisbis}
\end{equation}
\end{proof}

\begin{lem}
\label{lem:lap-conf}
    Let $u, w \in C^\infty_{>0}(N)$ and $r \in \mathbb{R}$. Then,
    \begin{equation}
     u^{r}\frac{1}{w}\Delta_{u^{r}h}w = r \frac{n-2}{4}
     \frac{\Delta_{h}(uw)}{uw} - r \frac{n-2}{4}
\frac{\Delta_{h}u}{u} + \left(1 - r \frac{n-2}{4}\right)
\frac{\Delta_{h}w}{w}.
        \label{eq:lap-conf1}
    \end{equation}
    In particular, if $w=u$, then
    \begin{equation}
     u^{r}\frac{1}{u}\Delta_{u^{r}h}u = r \frac{n-2}{4}
     \frac{\Delta_{h}u^{2}}{u^{2}} + \left(1 - r \frac{n-2}{2}\right)
\frac{\Delta_{h}u}{u} = \frac{1}{r\frac{n-2}{2}+1}
     \frac{\Delta u^{r\frac{n-2}{2}+1}}{u^{r\frac{n-2}{2}+1}},
        \label{eq:lap-conf2}
    \end{equation}
    where the latter equality is true when $\displaystyle r
    \neq -\frac{2}{n-2}$.
    Moreover, if $n \geq 3 $ and $\displaystyle r=\frac{4}{n-2}$, then
    \begin{equation}
     u^{\frac{4}{n-2}}\frac{1}{u}\Delta_{u^{\frac{4}{n-2}}h}u=
     \frac{\Delta_{h}u^{2}}{u^{2}} -
\frac{\Delta_{h}u}{u}=\frac{1}{3}\frac{\Delta_{h}u^{3}}{u^{3}}.
        \label{eq:lap-conf3}
    \end{equation}
\end{lem}

\begin{proof}
First of all, we observe that
\begin{equation}
     \frac{\Delta_{h}(uw)}{uw} = 2 h\left(\frac{\nabla
u}{u},\frac{\nabla w}{w}\right) +
     \frac{\Delta_{h}u}{u} + \frac{\Delta_{h}w}{w},
        \label{eq:lap-product1}
    \end{equation}
hence
\begin{equation}
     h\left(\frac{\nabla u}{u},\frac{\nabla
     w}{w}\right)=
     \frac{1}{2}\frac{\Delta_{h}(uw)}{uw}-
     \frac{1}{2}\left( \frac{\Delta_{h}u}{u} +
     \frac{\Delta_{h}w}{w}\right).
        \label{eq:lap-product2}
\end{equation}
On the other hand, by \textbf{Lemma
\ref{lem:conformal-laplacian-bis}},
\begin{equation}
 \begin{array}{lll}
  \displaystyle u^{r}\frac{1}{w}\Delta_{u^{r}h}w &=&\displaystyle
  r {\frac{n-2}{2}} \, h\left(\frac{\nabla u}{u},\frac{\nabla w}{w}
\right)
  + \frac{\Delta_{h}w}{w}\\
   &=& \displaystyle r {\frac{n-2}{4}}\left[\frac{\Delta_{h}(uw)}{uw}
- \frac{\Delta_{h}u}{u}
   -
     \frac{\Delta_{h}w}{w}\right] + \frac{\Delta_{h}w}{w}\\
   &=& \displaystyle r {\frac{n-2}{4}} \frac{\Delta_{h}(uw)}{uw} - r
{\frac{n-2}{4}}
   \frac{\Delta_{h}u}{u} +
   \left(1-r\frac{n-2}{4}\right)\frac{\Delta_{h}w}{w}.
  \end{array}
  \label{eq:lap-conf3-bis}
\end{equation}

In \eqref{eq:lap-conf2}, the first equality is immediate by taking
$w=u$ in \eqref{eq:lap-conf1}. In order to obtain the second
equality of \eqref{eq:lap-conf2} it is sufficient to apply Remark
\ref{rem:2.1} with $\displaystyle \alpha = \beta =
\frac{1}{r\frac{n-2}{2}+1} $. Finally, \eqref{eq:lap-conf3} is an
obvious consequence of \eqref{eq:lap-conf2}.
\end{proof}

\begin{rem} Now we compute the useful relation between the scalar
curvatures under a conformal change in the metric $h$
when the conformal metric is written in the form
$\tilde{h}=v^{r}h$, $h \in C^\infty_{>0}(N)$ instead of an
exponential form like in \eqref{eq:yamabe-1}. Consider $v^{r} =
e^{\eta}$, so that $\eta = r \log v$ and applying
\eqref{eq:yamabe-1} and \eqref{eq:lappotencia} (note that $t\neq
0, 1 $) we obtain
\begin{equation}\label{eq:conf 1}
  \begin{array}{lll}
     v^{r} S_{\tilde{h}} & = & \displaystyle S_{h} - (n-1)r
\left[\Delta_{h} \log v +
  \frac{n-2}{4} |\nabla \log v|^{2}\right]\\
     & = & \displaystyle S_{h} - (n-1)r \left[\left(-1 +
\frac{n-2}{4}r \right)\frac{|\nabla v|^{2}}{v^{2}} +
      \frac{\Delta_{ h} v}{v}\right]\\
     & = & \displaystyle S_{h} - (n-1)r \left[\left(-1 +
\frac{n-2}{4}r
     \right)\frac{1}{(t-1)t}
     \frac{\Delta_{h}v^{t}}{v^{t}}\right.\\
     & + & \left. \displaystyle
     \left(1-\left(-1 + \frac{n-2}{4}r
     \right)\frac{1}{t-1}\right)\frac{\Delta_{h}v}{v} \right].
  \end{array}
\end{equation}
Without lose of generality we assume $r$ is nonzero, it is clear
that $S_{\tilde{h}} = \displaystyle S_{h}$ when $r = 0$. At this
point, we have two cases:
\begin{description}
  \item[$(n \geq 3)$]
By Remark \ref{rem:2.1} with $\displaystyle \alpha = \beta =
\frac{4}{n-2} \frac{1}{r} $,
\begin{equation}\label{eq:conf 2}
  v^{r} S_{\tilde{h}} = \displaystyle S_{h} -
  (n-1)\frac{4}{n-2}\frac{\Delta_{h}v^{\frac{n-2}{4}r}}{
v^{\frac{n-2}{4}r}
  },
\end{equation}
which contents as a particular case \eqref{eq:yamabe} when
$\displaystyle r = \frac{4}{n-2}$.
  \item[$(n=2)$]
In this case \eqref{eq:conf 1} says
\begin{equation}\label{eq:conf 3}
  v^{r} S_{\tilde{h}} = \displaystyle S_{h} -
  \frac{r}{t-1}\left[-\frac{1}{t}\frac{\Delta_{h}v^{t}}{v^{t}} + t
\frac{\Delta_{h}v}{v}
  \right].
\end{equation}
\noindent Moreover, if we apply \eqref{eq:operator} the latter
equation becomes
\begin{equation}\label{eq:conf 4}
  v^{r} S_{\tilde{h}} = \displaystyle S_{h} +
  r \left(\frac{|\nabla v|^{2}}{v^{2}}-\frac{\Delta v}{v}\right).
\end{equation}
Note that in \eqref{eq:conf 3} it is not possible to apply Remark
\ref{rem:2.1}.
\end{description}
\end{rem}

\begin{rem} \label{rem: computations with conformal change} Now, as we
mentioned in \S 1, we will outline an alternative proof of
\textbf{of Theorem \ref{thm:scurv conf warped m ge 2}} by applying
a conformal change metric technique like in
\cite{DobarroLamiDozo87}. We will concentrate in \textbf{Theorem
\ref{thm:scalar sbcwp m ge 2}} when $m \ge 3$. The same order of
ideas may be used for the case $m=2$.
\end{rem}

%
%
%


\begin{proof}~\Big( of \textbf{Theorem
\ref{thm:scalar sbcwp m ge 2}} when $m \ge 3$, $ \displaystyle \mu
\neq -\frac{1}{m-2}$ \Big)
Since $g = \tilde{g_{B}}+\psi^{2}g_{F} $ with $\tilde{g_{B}}=
\psi^{2\mu}g_{B}$, an application of \eqref{eq:scalar curv wp-1}
to $\psi $ results
\begin{equation}
S=-2k \frac{\Delta_{\psi^{2\mu}g_{B}}\psi}{\psi} - k(k-1)\frac{
\psi^{2\mu} g_{B} (\psi^{-2\mu}\nabla \psi,\psi^{-2\mu}\nabla \psi
)
}{\psi^{2}} + S_{\psi^{2\mu} g_{B}} + \frac{S_{g_{F}}}{\psi^{2}}.
    \label{eq:scalar total m ge 3-1}
\end{equation}
So by multiplying by $\psi^{2\mu}$ and applying ~
\eqref{eq:lap-conf2} (note that $ \displaystyle \mu \neq
-\frac{1}{m-2}$), \eqref{eq:lappotencia} (with $t\neq 0, 1 $) and
equation \eqref{eq:conf 2} we obtain
\begin{equation}
  \begin{array}{lll}
\psi^{2\mu}S & = & \displaystyle -\frac{2k}{\mu (m-2) + 1}
\frac{\Delta_{g_{B}}\psi^{\mu (m-2)+1}}
{\psi^{\mu (m-2)+1}}\\
&&- k(k-1)  \displaystyle\frac{1}{t-1}\left[\frac{1}{t}
\frac{\Delta_{g_{B}}
\psi^{t} }{\psi^{t}} - \frac{\Delta_{g_{B}} \psi }{\psi}\right]\\
%
%
%
%
& \quad & +\displaystyle S_{g_{B}} -
  (m-1)\frac{4}{m-2}\frac{\Delta_{h}\psi^{\frac{m-2}{4}2\mu}}{
\psi^{\frac{m-2}{4}2\mu}
 }
+ S_{g_{F}}\psi^{2(\mu - 1)} \\
& = & -\displaystyle \left[\frac{2k}{\mu (m-2) + 1}
\frac{\Delta_{g_{B}}\psi^{\mu (m-2)+1}}
{\psi^{\mu (m-2)+1}} +   \frac{k(k-1)}{t(t-1)}
\frac{\Delta_{g_{B}} \psi^{t} }{\psi^{t}} \right.\\
&&\displaystyle -
\frac{k(k-1)}{(t-1)} \frac{\Delta_{g_{B}} \psi }{\psi} \\
& \quad &  \left. \displaystyle  +
  (m-1)\frac{4}{m-2}\frac{\Delta_{h}\psi^{\frac{m-2}{4}2\mu}}{
\psi^{\frac{m-2}{4}2\mu}
 }\right]
+ S_{g_{B}} + S_{g_{F}}\psi^{2(\mu - 1)}.
   \label{eq:scalar total 2}
  \end{array}
\end{equation}
The hypothesis of Lemma \ref{lem:2.1} is verified, indeed: since
$\displaystyle \mu\neq -\frac{k}{m-1}$,
\begin{equation}
  2(k+(m-1)\mu)\neq 0
\end{equation}
and
\begin{equation}
\label{eq:alpha condition}
  \begin{array}{rrr}
   2k (\mu (m-2) + 1) + k(k-1) + (m-1)(m-2)\mu^{2} &=& \\
  \{[k+(m-1)\mu]+(1-\mu)\}k + (m-2)\mu [k+(m-1)\mu] &=& \\
    \quad[k+(m-1)\mu](k+(m-2)\mu)+(1-\mu)k &=& \\
    (m-1)(m-2)\mu^{2} + 2(m-2)k\mu +(k+1)k &>& 0
  \end{array}
\end{equation}
Thus, by applying Lemma \ref{lem:2.1} with \\
$\displaystyle \alpha =
\frac{2[k+(m-1)\mu]}{\{[k+(m-1)\mu]+(1-\mu)\}k + (m-2)\mu
[k+(m-1)\mu]}$ and also \\ $\beta = \alpha 2[k+(m-1)\mu] $ (thus
$\beta
> 0$) and $u=\psi^{\frac{1}{\alpha}}$, we obtain that:

\begin{equation}
   -\beta \frac{\Delta_{g_{B}}u}{u}
   =
   u^{2 \mu \alpha} S - S_{g_{B}} - S_{g_{F}}u^{2(\mu -
   1)\alpha}.
        \label{eq:scalar total 3}
\end{equation}

\end{proof}

\begin{rem} By using the latter technique,
the case of $\displaystyle \mu = -\frac{k}{m-2}$ must be analyzed
separately. However, it is possible to prove \eqref{eq:scalar
curvature m ge 2}
in a similar way too.
\end{rem}

\providecommand{\bysame}{\leavevmode\hbox
to3em{\hrulefill}\thinspace}

\end{document}